\documentclass[reqno]{amsart}
\usepackage{amsfonts,amssymb,amsmath,fontenc}
\usepackage{latexsym,wasysym,mathrsfs,enumerate}
\usepackage{graphics}
\usepackage{color,graphicx}
\usepackage{hyperref}

\definecolor{darkred}{rgb}{.0,.0,.45}
\definecolor{darkgreen}{rgb}{.0,.55,0}
\definecolor{darkredd}{rgb}{.8,0,0}
\hypersetup{
 colorlinks,%
       citecolor=darkgreen,%
    filecolor=magenta,%
    linkcolor=darkredd,%
    urlcolor=darkred,%
}

\DeclareMathAlphabet{\mathpzc}{OT1}{pzc}{m}{it}

\begin{document}

\pagenumbering{arabic}
\title[Stability of Solitary Waves of the BO--ZK
Equation] {Stability and Decay properties of Solitary-wave solutions
to the generalized BO--ZK equation}



\maketitle

\begin{center}
 {\large \textbf{Amin Esfahani}}\\ {\small School of Mathematics and Computer Science \\
Damghan University\\
 Damghan, Postal Code 36716--41167, Iran \\
  E-mail: amin@impa.br, esfahani@du.ac.ir,  }\\
 \end{center}

\begin{center}
 {\large \textbf{Ademir Pastor}}\\ {\small IMECC--UNICAMP\\
Rua S\'ergio Buarque de Holanda, 651\\ Cidade Universit\'aria,
13083-859, Campinas--SP, Brazil.\\
 E-mail: apastor@ime.unicamp.br}\\
 \end{center}

 \begin{center}
 {\large \textbf{Jerry L. Bona}}\\ {\small Department of Mathematics, Statistics and Computer Science \\
University of Illinois at Chicago\\ 851 S. Morgan Street MC 249
Chicago, Illinois 60601, USA.\\
 E-mail: bona@math.uic.edu}
 \end{center}

\newtheorem{theorem}{{Theorem}}[section]
\newtheorem{lemma}[theorem]{{\sc \textbf{Lemma}}}
\newtheorem{definition}[theorem]{{\sc \textbf{Definition}}}
\newtheorem{proposition}{\begin{small}Proposition\end{small}}
\newtheorem{remark}[theorem]{{\sc \textbf{Remark}}}
\newtheorem{corollary}[theorem]{{\sc \textbf{Corollary}}}
\newcommand{\rr}{{\mathbb R}}
\newcommand{\dl}{\Delta}
\newcommand{\dd}{{\rm d}}
\newcommand{\T}{{\mathbb T}}
\newcommand{\C}{{\mathbb C}}
\newcommand{\N}{{\mathbb N}}
\newcommand{\lam}{{\lambda}}
\newcommand{\ff}{\varphi}
\newcommand{\ti}{\widetilde}
\newcommand{\what}{\widehat}
\newcommand{\Z}{{\mathbb Z}}
\newcommand{\om}{\omega}
\newcommand{\dr}{S_r(\xi)}
\newcommand{\dk}{S_k}
\newcommand{\sech}{\textmd{sech}}

\begin{abstract}
Studied here is the generalized Benjamin-Ono--Zakharov-Kuznetsov  equation
\begin{equation}
u_t+u^pu_x+\alpha\mathscr{H}u_{xx}+\varepsilon u_{xyy}=0, \quad
(x,y)\in\rr^2\!,\;\;t\in \rr^+\!,
\end{equation}
 in two space dimensions. Here, $\mathscr{H}$ is the Hilbert
transform and subscripts denote partial differentiation.    We classify when equation (1) possesses solitary-wave
solutions in terms of the signs of the constants $\alpha$ and $\varepsilon$ appearing in
 the dispersive terms and the strength of the
nonlinearity. Regularity and decay properties of these solitary wave are determined
and their stability is studied.
\end{abstract}

\numberwithin{equation}{section}

\section{Introduction}
This paper is concerned with existence and non-existence, stability
and some decay properties of solitary-wave solutions of the
two-dimensional generalized Benjamin-Ono--Zakharov-Kuznetsov equation
(BO--ZK equation henceforth),
\begin{equation}
u_t+u^pu_x+\alpha\mathscr{H}u_{xx}+\varepsilon u_{xyy}=0, \quad
(x,y)\in\rr^2\!,\;\;t\in\rr^+\!.\label{2bozk-main}
\end{equation}
Here $p>0, \alpha$ and $\varepsilon$  are non-zero real constants with $\varepsilon$ normalized to $\pm1$ by
appropriately rescaling the $y$-variable while $\mathscr{H}$ is the Hilbert transform
\[
\mathscr{H}u(x,y,t)=\mathrm{p.v.}\frac{1}{\pi}\int_\rr\dfrac{u(z,y,t)}{x-z}\;dz,
\]
in the
$x$-variable, where $\mathrm{p.v.}$ denotes the Cauchy principal value.

When $p=1$, this equation arises as a model for electromigration in thin
nanoconductors on a dielectric substrate (see \cite{jcms,lmsv}).   Equation
(\ref{2bozk-main}) may also be viewed as one of the natural, two-dimensional
generalizations of the one-dimensional Benjamin-Ono equation
 in much the same way that the Kadomtsev-Petviashvili equation
and the Zakharov-Kuznetsov equation generalize the Koreteweg-de Vries equation.

The generalized Benjamin-Ono equation
$$
u_t+u^pu_x+\alpha\mathscr{H}u_{xx}=0, \hspace{1.4cm}
x\in\rr,\;\;t\in\rr^+\!, \label{bo}
$$
and its counterpart
$$
u_t+u^pu_x+\alpha\mathscr{H}u_{xx} + \beta u_{xxx} =0, \quad
x\in\rr,\;\;t\in\rr^+\!, \label{ben}
$$
taking into account surface tension effects between the two layers
of fluid, have been considered by many authors.
Well-posedness issues for the pure
 initial-value problem have attracted a lot of interest recently (see, {\it e.g.}
\cite{burq-planchon,ken-kon,koch-tzvetkov,molrib,ponce,tao}).  Questions
about the existence and stability of solitary traveling-waves have been
investigated in \cite{Al}--\cite{BK} and \cite{KB}.

Theory for the generalized
Zakharov-Kuznetsov equation
$$
u_t+u^pu_x+\alpha u_{xxx}+\varepsilon u_{xyy}=0, \quad
(x,y)\in\rr^2,\;\;t\in\rr^+\!, \label{zk}
$$
is less abundant. Well-posedness was studied in \cite{fa, falipa, grhe,mopi,lp1,rv}.  As far as we know, the only results
concerning existence and nonlinear stability of solitary-wave
solutions of this equation was provided in \cite{db}.

The solitary-wave solutions of interest here have the form
$u(x,y,t)=\ff(x-ct,y)$, where $ c \neq 0$ is the speed of propagation
 and $u$ belongs to a natural
function space denoted $\mathscr{Z}$ and introduced presently.
Substituting this form into (\ref{2bozk-main}),
integrating once with respect to the variable $z = x-ct$ and assuming $\ff(z,y)$ decays suitably for large
values of $|z|$, it transpires that $\ff$ must satisfy
\begin{equation}\label{2bozk}
-c\ff+\frac{1}{p+1}\ff^{p+1}+\alpha\mathscr{H}\ff_x+\varepsilon
\ff_{yy}=0,
\end{equation}
where we have replaced the variable $z$ by $x$.

\begin{remark}   \label{scalingrem}
When it is convenient, it may be assumed that \eqref{2bozk} has the
norma\-lized form
\begin{equation}\label{3bozk}
-\ff+\frac{1}{p+1}\ff^{p+1}+\mathscr{H}\ff_x\pm \ff_{yy}=0
\end{equation}
by scaling the independent and dependent variables, \textit{viz.}
$$
u(x,y,t)=av(bx,dy,et)
$$
where $a^p=c$, $e=b=c/\alpha$ and $d=\varepsilon/c^2$. If instead,
we insist that $d>0$, so $\varepsilon=+1$, then equation
\eqref{2bozk} may be taken in the form
\begin{equation}\label{3bozk1}
-\ff+\frac{1}{p+1}\ff^{p+1}\pm\mathscr{H}\ff_x+ \ff_{yy}=0.
\end{equation}
Of course, throughout, it will be presumed that the power $p$ appearing
in the nonlinearity is rational and has the form $k/m$ where $k$ and
$m$ are relatively prime and $m$ is odd.  This restriction
allows us to define a  branch of the mapping $w \mapsto w^{\frac{1}{m}}$
that is real on the real axis.
\end{remark}

Attention is now turned to the structure of the paper. The theory
begins by examining when solitary-wave solutions of
\eqref{2bozk-main} exist. As pointed out in
\cite{lmsv}, no exact formulas are known for solitary-wave solutions
to (\ref{2bozk-main}), so an existence theory is needed before
questions of stability can be addressed.  Pohojaev-type identities are used to show
that solitary-wave solutions do not exist for certain values of $p$
and signs of $\varepsilon$ and $\alpha$. In some of the cases where
such solutions are not prohibited by elementary inequalities, a
suitable minimization problem can be solved using Lions'
concentration-compactness principle \cite{l1,l2} (see Theorem
\ref{nonexistence}). For example, our results imply there are
solitary-wave solutions when $c>0$, $\alpha<0$, $\varepsilon>0$ and
$0<p<4$. Moreover, these solutions are shown to be ground states.

With solitary waves in hand, their orbital stability is at issue.
The variational approach of Cazenave and Lions \cite{cl} comes to
the fore in Section \ref{stability-sec} in establishing stability
for the case $\alpha\varepsilon<0$, $c\alpha<0$, and $0<p<4/3$.
Complementary instability results  appeared in  \cite{prep-1} for
the same conditions on $c,\alpha$ and $\varepsilon$, but with
$4/3<p<4$.

The regularity and decay properties of the solitary-wave solutions
shown to exist in Section \ref{non(exist)} are developed in Sections
\ref{decay and regul} and \ref{decaysec}. Solitary-wave solutions are shown to be
positive and real analytic. They are symmetric about their peak with
respect to both the direction of propagation and the transverse
direction. Moreover, solitary waves decay   to zero algebraically in
the direction of propagation and exponentially  in the
transverse direction. Some of the results in Section \ref{decay and
regul} inform the analysis of instability in \cite{prep-1}.

In the theory developed here, the issue of well-posedness is not
addressed. The presumption throughout is that suitable
well-posedness obtains for these models. Detailed analysis of the
initial-value problem  appeared in \cite{CP} and \cite{prep-4}.

\begin{remark}
The scale-invariant Sobolev spaces for the BO--ZK equation
(\ref{2bozk-main}) are $\dot{H}^{s_1,s_2}(\rr^2)$, where
$2s_1+s_2=\frac{3}{2}-\frac{2}{p}$ (see the definitions below).
Hence a reasonable framework for studying local well-posedness of
the BO--ZK equation (\ref{2bozk-main}) is the family of spaces
${H}^{s_1,s_2}(\rr^2)$, $2s_1+s_2\geq\frac{3}{2}-\frac{2}{p}$.
\end{remark}

\begin{remark}
The $n$-dimensional version of \eqref{2bozk-main} is
\begin{equation}
u_t+u^pu_{x_1}+\alpha\mathscr{H}u_{x_1x_1}+\sum_{i=2}^n\varepsilon_i
u_{x_1x_ix_i}=0, \quad \label{nbozk}
\end{equation}
where $t\in \mathbb{R}^+$, $(x_1,x_2,\ldots,x_n) \in \mathbb{R}^n$
 and $\alpha, \varepsilon_i \in\mathbb{R}$, $i=2,\ldots,n$.
The theory developed here has natural analogs for \eqref{nbozk} which
will be developed later.
\end{remark}

\vskip.3cm

\noindent {\bf Notation and Preliminaries.} As already mentioned, the exponent $p$ in
\eqref{2bozk-main} is taken to be a rational number of the form $p =
k/m$, where $m$ and $k$ are relatively prime and $m$ is odd. This allows
the nonlinearity to be given a definition that is real-valued.  The notation $\what{f}=\what{f}(\xi,\eta)$ means the Fourier transform,
$$
\what{f}(\xi,\eta)=\int_{\rr^2}e^{-\mathrm{i}(x\xi+y\eta)}f(x,y)\;dxdy
$$
of $f=f(x,y)$.  For any $s\in \rr$, the space $H^s:=H^s(\rr^2)$ denotes the usual
isotropic, $L^2(\mathbb{R}^2)$-based, Sobolev space.
For $s_1,s_2\in\rr$, the anisotropic Sobolev space
$H^{s_1,s_2}:=H^{s_1,s_2}\left(\rr^2\right)$ is the set of all
distributions $f$ such that
\[
\|f\|_{H^{s_1,s_2}}^2=\int_{\rr^2}\left(1+
\xi^2\right)^{s_1}\left(1+\eta^2\right)^{s_2}|\widehat{f}
(\xi,\eta)|^2\;d\xi d\eta <\infty.
\]

\noindent The fractional Sobolev-Liouville spaces
$H_p^{(s_1,s_2)}:=H_p^{(s_1,s_2)}\left(\rr^2\right)$, $1\leq p
<\infty$, are the set of all functions $f \in L^p(\mathbb{R}^2)$
such that
\[
\|f\|_{H_p^{(s_1,s_2)}}=\|f\|_{L^p(\rr^2)}+
\sum_{i=1}^2\left\|D_{x_i}^{s_i}f\right\|_{L^p(\rr^2)}<\infty,
\]
where $D_{x_i}^{s_i}f$ denotes the Bessel derivative of order $s_i$
with respect to $x_i$ (see {\it e.g.} \cite{kolyada}, \cite{li}). For
short,  $H^{(k)}_p(\rr^2)$ denotes the space $H^{(k,k)}_p(\rr^2)$.

The particular space
$\mathscr{Z}:=H^{\frac12,0}\left(\rr^2\right)\cap
 H^{0,1}\left(\rr^2\right)=H^{\left(\frac12,1\right)}\left(\rr^2\right)$
arises naturally in the analysis to follow.  It can be characterized
alternatively as the closure of
 $C_0^\infty(\rr^2)$ with respect to the norm
\begin{equation}
\|\ff\|_\mathscr{Z}^2=\|\ff\|_{L^2(\rr^2)}^2+\left\|\ff_y\right\|_{L^2(\rr^2)}^2+
\left\|D_x^{1/2}\ff\right\|_{L^2(\rr^2)}^2,
\end{equation}
where $D_x^{1/2}\ff$ denotes the fractional derivative of order
$1/2$ with respect to $x$, defined via its Fourier transform by
$\widehat{D_x^{1/2}\ff}(\xi,\eta)=|\xi|^{1/2}\widehat{\ff}(\xi,\eta)$.

\begin{remark}\label{embdd-lem-1}
By combining  fractional Gagliardo-Nirenberg and H\"older's inequality one can deduce the existence of a positive constant $C$ such that
\begin{equation}\label{inhomo-embed}
\|u\|_{L^{p+2}}^{p+2}\leq C{\|u\|}_{L^2}^{(4-p)/2}{\|D_x^{1/2}
u\|}_{L^2}^{p}{\|u_y\|}_{L^2}^{p/2}, \quad 0\leq p<4.
\end{equation}
 This in turn implies the continuous embedding
 \begin{equation}\label{embedding-0}
\mathscr{Z}\hookrightarrow L^p\left(\rr^2\right), \quad 0\leq p<4.
\end{equation}
\end{remark}

\section{Solitary waves}\label{non(exist)}

This section is devoted to establishing  existence and non-existence results for
 solitary-wave solutions of the  BO-ZK equations. We begin with a non-existence result.

\begin{theorem}\label{nonexistence}
Equation (\ref{2bozk}) cannot have a non-trivial solitary-wave
solution unless either
\begin{enumerate}[(i)]
\item  $\varepsilon=1$, $c>0$, $\alpha<0$, $p<4$,
\item $\varepsilon=-1$, $c<0$, $\alpha>0$, $p<4$,
\item $\varepsilon=1$, $c<0$, $\alpha<0$, $p>4$, or
\item $\varepsilon=-1$, $c>0$, $\alpha>0$, $p>4$.
\end{enumerate}
\end{theorem}
\begin{proof}
 This follows from some Pohojaev-type identities.
If
(\ref{2bozk}) is multiplied by $\ff$, $x\ff_x$ and $y\ff_y$ and
the results integrated over $\rr^2$, then the identities
\begin{gather}
\int_{\rr^2}\left(-c\ff^2+\alpha\ff\mathscr{H}\ff_x-\varepsilon\ff_y^2
+\frac{1}{p+1}\ff^{p+2}\right) \! dxdy=0\label{2bozk-nonex-1},\\
\int_{\rr^2}\left(c\ff^2+\varepsilon\ff_y^2-\frac{2}{(p+1)(p+2)}
\ff^{p+2}\right)\!dxdy=0\label{2bozk-nonex-2},\\
\int_{\rr^2}\left(c\ff^2-\alpha\ff\mathscr{H}\ff_x-\varepsilon\ff_y^2-
\frac{2}{(p+1)(p+2)}\ff^{p+2}\right)\! dxdy=0,\label{2bozk-nonex-3}
\end{gather}
emerge.  These
formulas follow from the elementary
properties of the Hilbert transform together with suitably chosen formal  integrations by parts. The identities can be justified
for functions of the minimal regularity required for them to make
sense by first establishing them for smooth solutions and then using a
 standard truncation argument as in \cite{dbs1}.

Summing \eqref{2bozk-nonex-1} and \eqref{2bozk-nonex-2} leads to
\begin{equation}
\int_{\rr^2}\left(\alpha\ff\mathscr{H}\ff_x
+\frac{p}{(p+1)(p+2)}\ff^{p+2}\right)\!dxdy=0\label{2bozk-nonex-5},
\end{equation}
whilst adding \eqref{2bozk-nonex-2} and \eqref{2bozk-nonex-3} yields
\begin{equation}
\int_{\rr^2}\left(c\ff^2-\frac{\alpha}{2}\ff\mathscr{H}\ff_x
-\frac{2}{p+1}\ff^{p+2}\right)\! dxdy=0\label{2bozk-nonex-6}.
\end{equation}
If the integral of $\ff^{p+2}$ is eliminated between \eqref{2bozk-nonex-5} and
\eqref{2bozk-nonex-6}, there appears
\begin{equation}
\int_{\rr^2}\left(2pc\ff^2+\alpha(4-p)\ff\mathscr{H}\ff_x\right)\!
dxdy=0\label{2bozk-nonex-7}.
\end{equation}
On the other hand, adding \eqref{2bozk-nonex-1} and
\eqref{2bozk-nonex-3} gives
\begin{equation}
\int_{\rr^2}\left(2\varepsilon\ff_y^2
-\frac{p}{(p+1)(p+2)}\ff^{p+2}\right)\!dxdy=0\label{2bozk-nonex-8}.
\end{equation}
Finally, substituting  \eqref{2bozk-nonex-2} into
\eqref{2bozk-nonex-8}, there obtains
\begin{equation}
\int_{\rr^2}\left(pc\ff^2+\varepsilon(p-4)\ff_y^2\right)\!dxdy=0
\label{2bozk-nonex-9}.
\end{equation}
The advertised results  follow immediately from
(\ref{2bozk-nonex-7}) and (\ref{2bozk-nonex-9}).
 \end{proof}

For cases (i) and (ii) from Theorem \ref{nonexistence}, the existence
of solitary-wave solutions of \eqref{2bozk-main} is established
in the next result.

\begin{theorem}\label{existence}
Let $\alpha\varepsilon, c\alpha<0$ and $p=\frac{k}{m}<4$, where
$m\in\N$ is odd and $m$ and $k$ are relatively prime. Then equation
(\ref{2bozk}) admits a non-trivial solution $\ff\in \mathscr{Z}$.
\end{theorem}
 \begin{proof} The proof is based on the
concentration-compactness principle \cite{l1,l2}. Suppose that
$\alpha<0$ (the proof for $\alpha>0$ is similar). Without loss of
generality, assume that $\alpha=-1$ and $c=1$ so that
$\varepsilon=+1$ (see Remark \ref{scalingrem}) and consider the
minimization problem
\begin{equation} \label{minimiz-constr-lp}
I_\lam=\inf\left\{ I(\ff) \;;\;\ff\in\mathscr{Z}\;,\;J(\ff)=
\int_{\mathbb{R}^2}\ff^{p+2}dxdy=\lam\right\}
\end{equation}
where $\lam\neq0$ and
\[
I(\ff)=\frac{1}{2}\int_{\rr^2}\left(\ff^2+\ff\mathscr{H}\ff_x+\ff_y^2\right)dxdy=\frac{1}{2}\|\ff\|_{\mathscr{Z}}^2.
\]
Clearly, $I_\lam<\infty$ if there are elements $\ff\in\mathscr{Z}$ such that
$\int_{\mathbb{R}^2}\ff^{p+2}\;dxdy=\lam. $\begin{footnote}{Depending on $p$, this might
require that $\lambda > 0$.  Of course, $I_\lambda$ is a number, but we will sometimes
refer to it as the minimization problem. For example, the phrase ``$\{\phi_n\}$ is a
minimizing sequence for the problem $I_\lambda$" means that $J(\phi_n) = \lambda$ for all $n$
and  $I(\phi_n) \to I_\lambda$ as $n \to \infty$.}\end{footnote}  The embedding
\eqref{embedding-0} allows us to adduce a positive constant $C$ such
that
\[
0<|\lam|=\left|\int_{\rr^2}\ff^{p+2}\;dxdy\right|\leq
C\|\ff\|_\mathscr{Z}^{p+2}= CI(\ff)^{\frac{p+2}{2}},
\]
from which one concludes that $I_\lambda\geq
\left(\frac{|\lam|}{C}\right)^{\frac{2}{p+2}}>0$.

For suitable $\lambda$ let  $\{\ff_n\}_{n \in \N}$ be a minimizing sequence for $I_\lam$. For $n=1,2,\cdots$ and  $r>0$,  define the
concentration function $Q_n(r)$ associated to $\ff_n$ by
\[
Q_n(r)=\sup_{(\ti{x},\ti{y})\in\rr^2}\int\limits_{B_r(\ti{x},\ti{y})}\rho_n\,dxdy
\]
where
$\rho_n=\left|\ff_n\right|^2+\left|D_x^{1/2}\ff_n\right|^2+\left|\partial_y\ff_n\right|^2$
and $B_r(x,y)$ denotes the ball of radius $r>0$ centered at
$(x,y)\in\rr^2$. If  evanescence of the sequence $\{\ff_n\}_{n \in \N}$ occurs,
which is to say, for any $r>0$,
\[
\lim_{n\to+\infty}\sup_{(\ti{x},\ti{y})\in\rr^2}\int\limits_{B_r(\ti{x},\ti{y})}\rho_n\,dxdy=0,
\]
then embedding \eqref{embedding-0} implies that
$\lim_{n\to\infty}\|\ff_n\|_{L^{p+2}}=0$, which contradicts the
constraint imposed for the minimization problem. Thus, according to
the concentration-compactness theorem, either dichotomy or
compactness must occur for the sequence $\{\ff_n\}_{n \in \N}$.

The occurrence of dichotomy is ruled out next. Suppose that $\gamma\in(0,I_\lam)$,
where it is assumed that
\[
\gamma=\lim_{r\to+\infty}\lim_{n\to+\infty}\sup_{(\ti{x},\ti{y})
\in\rr^2}\int\limits_{B_r(\ti{x},\ti{y})}\rho_n\,dxdy.
\]
By the definition of $\gamma$, for a given $\epsilon>0$, there exist
$r_1\in\rr$ and $N\in\N$ such that
\[
\gamma-\epsilon<Q_n(r)\leq Q_n(2r)<\gamma+\epsilon,
\]
for any $r\geq r_1$ and $n\geq N$. Hence, there is a sequence
$\{(\ti{x}_n,\ti{y}_n)\}_{n \in \N} \subset\rr^2$ for which
\[
\int\limits_{B_r(\ti{x}_n,\ti{y}_n)}\rho_n\;dxdy>\gamma-
\epsilon\quad \;\:
\mbox{and}\quad\:\: \int\limits_{B_{2r}(\ti{x}_n,\ti{y}_n)}\rho_n\;dxdy<\gamma+\epsilon.
\]
Let $\phi,\psi$ lie in  $C^\infty(\rr^2)$ and suppose
\begin{itemize}
\item $\mbox{supp}\;\phi\subset B_2(0,0)$,\;$\phi\equiv1$ on $B_1(0,0)$ and $0\leq\phi\leq1$,
\item $\mbox{supp}\;\psi\subset\rr^2\setminus B_1(0,0)$,\;$\psi\equiv1$ on $\rr^2\setminus B_2(0,0)$ and $0\leq\psi\leq1$.
\end{itemize}
Define the sequences $\{g_n\}_{n \in \N}$ and $\{h_n\}_{n \in \N}$ by
\[
g_n(x,y)=\phi_r((x,y)-(\ti{x}_n,\ti{y}_n))\ff_n\;\;\;\;\mbox{and}\;\;\;\; h_n(x,y)=\psi_r((x,y)-(\ti{x}_n,\ti{y}_n))\ff_n,
\]
where
\[
\phi_r(x,y)=\phi\left(\frac{(x,y)}{r}\right)\;\;\;\;\mbox{and}\;\;\;\; \psi_r(x,y)=\psi\left(\frac{(x,y)}{r}\right).
\]
It is clear that $g_n,h_n\in\mathscr{Z}$\!.

The following commutator estimate is helpful in obtaining the
splitting lemma to follow.

\begin{lemma}[\cite{cal,coffmey}]\label{coffmey}
Let $g\in C^\infty(\rr)$ with $g'\in L^\infty(\rr)$. Then
$[\mathscr{H},g]\partial_x$ is a bounded linear operator from $L^2(\rr)$ into $L^2(\rr)$ with
\[
\left\|[\mathscr{H},g]\partial_xf\right\|_{L^2(\rr)}\leq C\|g'\|_{L^\infty(\rr)}\|f\|_{L^2(\rr)}.
\]
\end{lemma}

The splitting lemma proved next enables us to rule out the
possibility of dichotomy occuring in the present context.

\begin{lemma}
Let $\{g_n\}_{n\in\N}$ and $\{h_n\}_{n\in\N}$ be as just defined.
Then, for every $\epsilon>0$, there exists
$\delta=\delta(\epsilon)>0$ with
$\lim_{\epsilon\to0}\delta(\epsilon)=0$, $\mu\in(0,I_\lam)$,
$n_0\in\N$ and $\rho\in(0,\lam)$ such that for all $n\geq n_0$,
\begin{gather}
|I(\ff_n)-I(g_n)-I(h_n)|\leq\delta,\label{split-1}\\
|I(g_n)-\mu|\leq\delta,\quad
|I(h_n)-I_\lam+\mu|\leq\delta,\label{split-2}\\
|J(\ff_n)-J(g_n)-J(h_n)|\leq\delta,\label{split-3}\\
|J(g_n)-\rho|\leq\delta,\quad|J(h_n)-\lam+\rho|
\leq\delta.\label{split-4}
\end{gather}
\end{lemma}
 \begin{proof}
 Obviously,  $\mathrm{supp}\;g_n \cap
\mathrm{supp}\; h_n = \emptyset$. Write
$g_n=\phi_r\ff_n$ and $h_n=\psi_r\ff_n$ so that
\[\begin{split}
2I(g_n)&=\int_{\rr^2}\phi_r^2\left[\ff_n^2+\ff_n\partial_x\mathscr{H}
\ff_n+\left(\partial_y^2\ff_n\right)^2\right]dxdy
+2\int_{\rr^2}\phi_r\ff_n(\partial_y\phi_r)(\partial_y\ff_n)dxdy\\
&\;\;\;+\int_{\rr^2}
\left[(\partial_y\phi_r)^2\ff_n^2+\ff_n\phi_r\mathscr{H}(\ff_n\partial_x\phi_r)
\right]dxdy+\int_{\rr^2}\ff_n\phi_r[\mathscr{H},\phi_r]\partial_x\ff_ndxdy
\end{split}\]
and
\[\begin{split}
2I(h_n)&=\int_{\rr^2}\psi_r^2\left[\ff_n^2+\ff_n\partial_x\mathscr{H}
\ff_n+\left(\partial_y^2\ff_n\right)^2\right]dxdy
+2\int_{\rr^2}\psi_r\ff_n(\partial_y\psi_r)(\partial_y\ff_n)dxdy\\
&\;\;\;+\int_{\rr^2}
\left[(\partial_y\psi_r)^2\ff_n^2+\ff_n\psi_r\mathscr{H}(\ff_n\partial_x\psi_r)\right]
dxdy+\int_{\rr^2}\ff_n\psi_r[\mathscr{H},\psi_r]\partial_x\ff_ndxdy.
\end{split}\]
Since $\|\phi_r\|_{L^\infty}=\|\psi_r\|_{L^\infty}=1$,
$\|\nabla\phi_r\|_{L^\infty}\leq\frac{1}{r}\|\nabla\phi\|_{L^\infty}$
and $\|\nabla\psi_r\|_{L^\infty}\leq\frac{1}{r}\|\nabla\psi\|_{L^\infty}$,
it follows from Lemma \ref{coffmey} that
\[
\left|I(g_n)-\frac{1}{2}\int_{\rr^2}\phi_r^2\left[\ff_n^2+\ff_n\partial_x\mathscr{H}
\ff_n+\left(\partial_y^2\ff_n\right)^2\right]\;dxdy\right|\leq\frac{1}{2}\delta(\epsilon)
\]
and
\[
\left|I(h_n)-\frac{1}{2}\int_{\rr^2}\psi_r^2\left[\ff_n^2+\ff_n
\partial_x\mathscr{H}\ff_n+\left(\partial_y^2\ff_n\right)^2\right]\;
dxdy\right|\leq\frac{1}{2}\delta(\epsilon).
\]
These inequalities imply (\ref{split-1}), from which, one infers
(taking subsequences if necessary) that there exists
$\mu=\mu(\epsilon)\in[0,I_\lam]$  such that
$\lim_{n\to\infty}I(g_n)=\mu$. In consequence, we see that
\[
|I(g_n)-I_\lam+\mu|\leq\delta(\epsilon).
\]
From (\ref{split-1}) again, the fact that $\mathrm{supp}\; g_n \cap
\mathrm{supp}\; h_n = \emptyset$ and  the embedding
(\ref{embedding-0}), one  obtains
\[
|J(\ff_n)-J(g_n)-J(h_n)|\leq C\delta(\varepsilon)
\]
for some constant $C$. It may therefore be presumed that there is a
$\rho=\rho(\epsilon)$ and
$\widetilde{\rho}=\widetilde{\rho}(\epsilon)$ such that
\[
\lim_{n\to+\infty}J(g_n)=\rho(\epsilon),\quad\lim_{n\to+\infty}J(h_n)=\ti{\rho}(\epsilon)
\]
with
$|\lam-\rho(\epsilon)-\ti{\rho}(\epsilon)|\leq\delta(\epsilon)$. If
$\lim_{\epsilon\to0}\rho(\epsilon)=0$, then for $\epsilon$
sufficiently small, it must be that  $J(h_n)>0$ for $n$ large
enough. Hence, by considering
$\left(\ti{\rho}(\epsilon)J(h_n)\right)^{\frac{1}{p+2}}h_n$, and
noting that
$J\left(\left(\ti{\rho}(\epsilon)J(h_n)\right)^{\frac{1}{p+2}}h_n\right)=\ti{\rho}(\epsilon)$,
it transpires that
\[
I_{\ti{\rho}(\epsilon)}\leq\liminf_{n\to+\infty}I(h_n)\leq I_\lam-\gamma+\delta(\epsilon),
\]
which leads to a contradiction since
$\lim_{\epsilon\to0}\ti{\rho}(\epsilon)=\lam$. Thus
$\rho=\lim_{\epsilon\to0}\rho(\epsilon)>0$. Necessarily $\rho<\lam$,
because the case $\rho=\lam$ is ruled out in the same manner as just
used to rule out $\rho=0$, but with $h_n$ replacing $g_n$ in the
argument. Since $\rho\in(0,\lam)$, one infers that necessarily
$\mu=\lim_{\epsilon\to+\infty}\mu(\epsilon)\in(0,I_\lam)$. This
completes the proof of the lemma.
 \end{proof}

Now, attention is returned to the proof that dichotomy cannot
happen. The previous lemma implies that
\begin{equation}
I_\lam\geq I_\rho+I_{\lam-\rho},
\end{equation}
which contradicts the subadditivity of $I_\lam$ coming from the fact
that $I_\lam=\lam^{2/(p+2)}I_1$. Hence dichotomy is ruled out.
\vspace{.2cm}

The remaining case in the concentration-compactness principle is
local compactness. Thus, there exists a sequence
$\{(x_n,y_n)\}_{n\in\N}\subset\rr^2$ such that for all $\epsilon>0$,
there are  finite values $R>0$ and $n_0>0$ with
\[
\int_{B_R(x_n,y_n)}\rho_n\,dxdy\geq\iota_\lam-\epsilon,
\]
for all $n\geq n_0$, where
\[
\iota_\lam=\lim_{n\to+\infty}\int_{\rr^2}\rho_n\,dxdy.
\]
This implies that for $n$ large enough,
\[
\int_{B_R(x_n,y_n)}|\ff_n|^2dxdy\geq\int_{\rr^2}|\ff_n|^2dxdy-2\epsilon.
\]
Since $\ff_n$ is bounded in the Hilbert space $\mathscr{Z}$, there
exists $\ff\in\mathscr{Z}$ such that a subsequence of
$\{\ff_n(\cdot-(x_n,y_n))\}_{n\in\N}$ (denoted again by
$\{\ff_n(\cdot-(x_n,y_n))\}_{n\in\N}$) converges weakly in
$\mathscr{Z}$ to $\ff$. It follows that
\[
\begin{split}
\int_{\rr^2}|\ff|^2\;dxdy&\leq\liminf_{n\to+\infty}\int_{\rr^2}|
\ff_n|^2\;dxdy \\
&\leq\liminf_{n\to+\infty}\int_{B_R(x_n,y_n)}|\ff_n|^2\;dxdy+2\epsilon\\
&=\liminf_{n\to+\infty}\int_{B_R(0,0)}|\ff_n((x,y)-(x_n,y_n))|^2\;dxdy+2\epsilon.
\end{split}
\]
But, when restricted to bounded sets in $\mathbb{R}^2$,
$\mathscr{Z}$ is compactly embedded into $L^2$. Consequently,
$\{\ff_n(\cdot-(x_n,y_n))\}_{n\in\N}$ may be presumed to converges strongly in the
Fr\'echet space $L^2_{loc}(\rr^2)$. The last inequality above
implies that this strong convergence also takes place in
$L^2(\rr^2)$ by what are, by now, standard arguments. Thus, because
of the embedding \eqref{embedding-0},
$\{\ff_n(\cdot-(x_n,y_n))\}_{n\in\N}$ also converges to $\ff$
strongly in $L^{p+2}(\rr^2)$, whence $J(\ff)=\lam$ and
\[
I_\lam=\lim_{n\to+\infty}I(\ff_n)=I(\ff),
\]
which is to say, $\ff$ is a solution of $I_\lam$.
\vspace{.2cm}

The Lagrange multiplier theorem now implies  there exists
$\theta\in\rr$ such that
\begin{equation}
\ff+\mathscr{H}\ff_x-\ff_{yy}=\theta(p+2)\ff^{p+1}
\end{equation}
as an equation in $\mathscr{Z}'$ (the dual space of $\mathscr{Z}$ in
$L^2-$duality). A change of scale yields a $\widetilde{\ff}$ which
satisfies (\ref{2bozk}).
 \end{proof}

\begin{remark}
Theorem \ref{existence} shows the existence of
solitary-wave solutions  of \eqref{2bozk-main} in the cases
$\mathrm{(i)}$ and $\mathrm{(ii)}$  in Theorem \ref{nonexistence}.
The question of existence or nonexistence of solitary waves in cases
$\mathrm{(iii)}$  and $\mathrm{(iv)}$ is currently open.
\end{remark}

\begin{definition}
A solution $\ff$ of equation (\ref{2bozk}) is called a ground state,
if $\ff$ minimizes the action
\[
\mathcal{S}(u)=\mathscr{E}(u)+c\mathscr{F}(u)
\]
among all solutions of (\ref{2bozk}), where
\[
\mathscr{F}(u)=\dfrac{1}{2}\int_{\rr^2}u^2\;dxdy
 \]
 and
\[
\mathscr{E}(u)=\dfrac{1}{2}\int_{\rr^2}\left(\varepsilon u_y^2- \alpha
u\mathscr{H}u_x-\dfrac{2}{(p+1)(p+2)}\;u^{p+2}\right)\;dxdy.
\]
\end{definition}

Next, it is established that the minima obtained in Theorem
\ref{existence} are precisely the ground-state solutions of
(\ref{2bozk}). The proof is inspired by that of Lemma 2.1 in
\cite{dbs2}.

\begin{theorem}\label{varitional-charac}
In the context of   equation  \eqref{2bozk} for solitary-wave solutions of the
BO-ZK equation, let
\[
\mathcal{K}(u)=\dfrac{1}{2}\int_{\rr^2}(cu^2+u_y^2)dxdy-\dfrac{1}{(p+1)(p+2)}J(u)
\]
with $J(u)=\int u^{p+2}dxdy$ as in \eqref{minimiz-constr-lp}.
Up to a change of scale, the following assertions about a function
$u^\ast\in\mathscr{Z}$ are equivalent:
\begin{enumerate}[(i)]
\item  If $J(u^\ast)=\lam^\ast$ then $u^\ast$ is a
minimizer of $I_{\lam^\ast}$,
\item  $\mathcal{K}(u^\ast)=0$ and
$$\inf\left\{\int_{\rr^2}u\mathscr{H}u_x\;dxdy,\;u\in\mathscr{Z},
\;u\neq0,\;\mathcal{K}(u)=0\right\}=\int_{\rr^2}u^\ast\mathscr{H}u^\ast_x\;dxdy,$$
\item  $u^\ast$ is a ground state,
\item  $\mathcal{K}(u^\ast)=0$ and
$$\inf\left\{\mathcal{K}(u),\;u\in\mathscr{Z},\;u\neq0,\;
\int_{\rr^2}u\mathscr{H}u_x\;dxdy=\int_{\rr^2}u^\ast\mathscr{H}u^\ast_x\;dxdy\right\}=0.$$
\end{enumerate}
\end{theorem}
\noindent\emph{Proof}.
We set $\lam^\ast=\left(2(p+1)I_1\right)^\frac{p+2}{p}$ and proceed with the proof.
\vspace{.2cm}

\indent
$\rm{(i)}\Rrightarrow\rm{(ii)}\!:$ \;Assume that $u^\ast$
satisfies \rm{(i)}.
Let $u\in\mathscr{Z}$ with $u\neq0$ and $\mathcal{K}(u)=0$,
from which it follows that $J(u)>0$.  Define
 \[u_\mu(x,y)=u\left(\dfrac{x}{\mu},y\right),\;\;\mbox{with}
 \;\;\;\mu=\dfrac{J(u^\ast)}{J(u)},\]
so that $J(u_\mu)=J(u^\ast)$ and $\mathcal{K}(u_\mu)=0.$
Since $u^\ast$ is a minimum of $I_{\lam^\ast}$, it must be the case that
$\mathcal{K}(u^\ast)=0$ and
\[
\mathcal{K}(u^\ast)+C_pJ(u^\ast)+\dfrac{1}{2}\int_{\rr^2}u^\ast\mathscr{H}u^\ast_x\;dxdy\leq
\mathcal{K}(u_\mu)+C_pJ(u_\mu)+\dfrac{1}{2}\int_{\rr^2}u_\mu\mathscr{H}(u_\mu)_x\;dxdy,
\]
where $C_p=\frac{1}{(p+1)(p+2)}$.  This in turn implies that
\[
\int_{\rr^2}u^\ast\mathscr{H}u^\ast_x\;dxdy\leq\int_{\rr^2}u\mathscr{H}u_x\;dxdy,
\]
and \rm{(ii)} holds.
\vspace{.2cm}

$\rm{(ii)}\Rrightarrow\rm{(iii)}:\;\;$ If $u^\ast$ satisfies
 \rm{(ii)}, then there is a Lagrange multiplier $\theta$ such that
\[
cu^\ast-u^\ast_{yy}+\theta\mathscr{H}u_x^\ast-\dfrac{1}{p+1}(u^\ast)^{p+1}=0.
\]
By multiplying the above equation by $u^\ast$\!, integrating
 by parts and using that $\mathcal{K}(u^\ast)=0$, we can see that $\theta$
is positive. Hence the scale change $u_\ast(x,y)=u^\ast(x/\theta,y)$
satisfies  equation (\ref{2bozk}).

On the other hand, the identity
$S(u)=\mathcal{K}(u)+\frac{1}{2}\int_{\rr^2}u\mathscr{H}u_x dxdy$
shows that if $u$ is a solution of (\ref{2bozk}), then
\[
S(u)=\dfrac{1}{2}\int_{\rr^2}u\mathscr{H}u_x\,dxdy\geq\dfrac{1}{2}\int_{\rr^2}u^\ast\mathscr{H}u^\ast_x\,dxdy
=\frac{1}{2}\int_{\rr^2}u_\ast\mathscr{H}(u_\ast)_x\,dxdy=S(u^\ast),
\]
whence $u^\ast$ is a ground state.
\vspace{.2cm}

$\rm{(iii)}\Rrightarrow\rm{(i)}:\;$ From the proof of Theorem
\ref{nonexistence}, one   sees that if $u$ is a solution of
(\ref{2bozk}), then $\mathcal{K}(u)=0$ and
\begin{equation}
I(u)=\dfrac{1}{2}\left(1+\dfrac{2}{p}\right)\int_{\rr^2}u\mathscr{H}u_x\,dxdy.
\end{equation}
Hence if $u^\ast$ is a ground state, then $u^\ast$ minimizes
 both $I(u)$ and $\int_{\rr^2}u\mathscr{H}u_xdxdy$ among all solutions of
(\ref{2bozk}). Let $\lam=J(u)$ and $\ti{u}$ be a minimum of $I_\lam$. Then
\begin{equation} \label{eqn}
I_\lam=I(\ti{u})\leq I(u^\ast)
\end{equation}
and there is a positive number $\theta$ such that
\[
c\ti{u}-\ti{u}_{yy}+\mathscr{H}\ti{u}_x=\dfrac{\theta}{p+1}\ti{u}^{p+1}.
\]
Using the equations satisfied by $\ti{u}$ and $u^\ast$\!,   inequality \eqref{eqn} is written as
\[
I_\lam=\dfrac{\lam\theta}{p+1}\leq\dfrac{\lam}{p+1},
\]
from which it is deduced immediately that $\theta\leq1$. On the other hand, $u_*=\theta^p\ti{u}$
satisfies   equation (\ref{2bozk}), and since $u^\ast$ is a ground
state, it must be the case that
\[
I(u^\ast)\leq I(u_*)\leq\theta^{2p}I(\ti{u}),
\]
so that $\theta\geq1$. In consequence,  $u^\ast=\ti{u}$ is a minimum
of $I_\lam$ with $\lam=\lam^\ast$.
\vspace{.2cm}

$\rm{(ii)}\Rrightarrow\rm{(iv)}:\;$ Let $u\in\mathscr{Z}$ with
$\int_{\rr^2}u\mathscr{H}u_xdxdy=\int_{\rr^2}u^\ast
 \mathscr{H}u^\ast_xdxdy$.
Suppose that $\mathcal{K}(u)<0$. Since $\mathcal{K}(\tau u)>0$ for
$\tau>0$ sufficiently small, then there is a $\tau_0\in(0,1)$ such
that $\mathcal{K}(\tau_0u)=0$. Thus by setting $\ti{u}=\tau_0u$, one
has $\ti{u}\in\mathscr{Z}$, $\mathcal{K}(\ti{u})=0$ and
\[
\int_{\rr^2}\ti{u}\mathscr{H}\ti{u}_x\,dxdy<\int_{\rr^2}u
\mathscr{H}u_x\,dxdy=\int_{\rr^2}u^\ast\mathscr{H}u^\ast_x\,dxdy,
\]
which contradicts \rm{(ii) and shows that $u^\ast$ satisfies
\rm{(iv) because $\mathcal{K}(u^\ast)=0$.
\vspace{.2cm}

\indent
$\rm{(iv)}\Rrightarrow\rm{(ii)}:\;$ Let $u\in\mathscr{Z}$ with
$\mathcal{K}(u)=0$ and $u\neq0$. Suppose that
\[
\int_{\rr^2}u\mathscr{H}u_x\;dxdy<\int_{\rr^2}u^\ast\mathscr{H}u^\ast_x\;dxdy.
\]
Since $\mathcal{K}(\tau u)<0$ for $\tau>1$, there is a
$\tau_0>1$ with $$\int_{\rr^2}(\tau_0u)\mathscr{H}
(\tau_0u)_x\,dxdy=\int_{\rr^2}u^\ast\mathscr{H}u^\ast_x\,dxdy$$
 and $\mathcal{K}(\tau_0u)<0.$
 This contradicts  \rm{(iv). Hence $\int_{\rr^2}u
 \mathscr{H}u_x\;dxdy\geq\int_{\rr^2}u^\ast\mathscr{H}u^\ast_x\;dxdy$ and \rm{(ii) holds.
 \hfill$\Box$\\

\begin{remark}
Note that the proof of the above theorem shows that, indeed,
$\rm{(i)}$ and $\rm{(iii)}$ are equivalent and imply $\rm{(ii)}$ and $\rm{(iv)}$,
which are also equivalent.   The converse holds modulo a scale change.
\end{remark}

\section{Stability}\label{stability-sec}

The notion of orbital stability employed here is the standard one.

\begin{definition}
Let $\ff_c$  be  a solitary-wave solution of (\ref{2bozk-main}). We
say that $\ff_c$ is orbitally stable if for all $\eta>0$, there is a
$\delta>0$ such that for any $u_0\in H^s\left(\rr^2\right)$, $s>2$,
with $\|u_0-\ff_c\|_\mathscr{Z}<\delta$, the corresponding solution
$u(t)$ of (\ref{2bozk-main}) with $u(0)=u_0$ satisfies
\[
\sup_{t\geq0}\inf_{r\in\rr^2}\|u(t)-\ff_c(\cdot-r)\|_\mathscr{Z}<\eta.
\]
\end{definition}

Some of the arguments below can be found in \cite{ang} where the stability of solitary waves for the generalized BO equation has been established.
Hereafter, without loss of generality, we take $\alpha=-1$ so
that $\varepsilon=+1$, and $c>0$.

The following theorem is a consequence of Theorem \ref{existence}
and it will be used to obtain the  stability results.

\begin{theorem}\label{stab-theo1}
Let $\lam\neq0$.
\begin{enumerate}[(i)]
\item Every minimizing sequence for the problem $I_\lam$ converges, up
to translations, in $\mathscr{Z}$ to an element in the set
 $$
M_\lam=\{\ff\in\mathscr{Z};\;I(\ff)=I_\lam,\; J(\ff)=\lam\}
$$
of minimizers for $I_\lam$.
\item  Let $\{\ff_n\}$ be a minimizing sequence for $I_\lam$. Then, it must be the case
that
\begin{gather}
\lim_{n\to+\infty}\inf_{\psi\in M_\lam,\;z\in\rr^2}\|\ff_n(\cdot+z)-\psi\|_{\mathscr{Z}}=0,\label{stab-eq-theo-1}\\
\lim_{n\to+\infty}\inf_{\psi\in M_\lam}\|\ff_n-\psi\|_{\mathscr{Z}}=0.\label{stab-eq-theo-2}
\end{gather}
\end{enumerate}
\end{theorem}
 \begin{proof} Part (i) follows immediately from the proof
of Theorem \ref{existence}. The equality \eqref{stab-eq-theo-1} is
proved by contradiction. Indeed, if \eqref{stab-eq-theo-1} does not
hold, then there exists a subsequence of the sequence $\{\ff_n\}$, denoted again
by $\{\ff_n\}$, and an $\epsilon>0$ such that
\[
\varpi=\inf_{\psi\in M_\lam, r\in\rr^2}\|\ff_n(\cdot+r)-\psi\|_\mathscr{Z}\geq\epsilon,
\]
for all $n$ sufficiently large. On the other hand, since $\{\ff_n\}$ is a minimizing
sequence for $I_\lam$, part  \rm{(i)} implies that there exists
a sequence $\{r_n\}\subset\rr^2$ such that, up to a subsequence,
$\ff_n(\cdot+r_n)\to\ff$ in $\mathscr{Z}$, as $n\to\infty$. Hence,
for $n$ large enough, it is inferred that
\[
\frac{\epsilon}{2}\geq\|\ff_n(\cdot+r_n)-\ff\|_\mathscr{Z}\geq\varpi\geq\epsilon,
\]
which is a contradiction.

The proof of (\ref{stab-eq-theo-2}) follows from
\eqref{stab-eq-theo-1}, the fact that if $\psi\in M_\lam$ then
$\psi(\cdot+r)\in M_\lam$ for all $r\in\rr^2$, and the
equality
\[\begin{split}
\inf_{\psi\in M_\lam}\|\ff_n-\psi\|_\mathscr{Z}=\inf_{\psi\in
M_\lam,r\in\rr^2} \|\ff_n-\psi(\cdot-r)\|_\mathscr{Z} =\inf_{\psi\in
M_\lam,r\in\rr^2}\|\ff_n(\cdot+r)-\psi\|_\mathscr{Z}.
\end{split}\]
This completes the proof of the theorem.  \end{proof}

The next lemma shows that there exists a $\lam>0$ such that
every element in the set of minimizers satisfies (\ref{2bozk}).

\begin{lemma}
If $\lam=\big(2(p+1)I_1\big)^\frac{p+2}{p}$ in the minimization
problem \eqref{minimiz-constr-lp}, then any $\ff\in M_\lam$   is
a solitary-wave solution of  (\ref{2bozk}).
\end{lemma}

For $\lam$ as in the preceding  lemma, define the set
\[
\mathscr{N}_c=\left\{\ff\in\mathscr{Z};\;J(\ff)=2(p+1)I(\ff)=\lam\right\}.
\]
It is clear that $M_\lam=\mathscr{N}_c$; the latter notation simply emphasizes
the dependence upon the wave speed $c$.   Next, for any $c>0$ and any
$\ff\in\mathscr{N}_c$,  define the function $d:\rr \to \rr$ by
\begin{equation} \label{functiond}
d(c)=E(\ff)+c\mathscr{F}(\ff).
\end{equation}

\begin{lemma}\label{stabl-lem-2}
The function $d$ in \eqref{functiond} is constant on
$\mathscr{N}_c$ and differentiable and strictly increasing for $c>0$.
Moreover, $d''(c)>0$ if and only if $0<p<\frac{4}{3}$.
\end{lemma}
 \begin{proof}
  It is straightforward to check that
\begin{equation} \label{a1}
d(c)=I(\ff)-\dfrac{1}{(p+1)(p+2)}J(\ff)=\dfrac{p}{2(p+1)(p+2)}J(\ff)=\dfrac{p(2(p+1))^\frac{2}{p}}{p+2}I_1^{\frac{p+2}{p}}\!.
\end{equation}
It is plain that $d$ is constant on $\mathscr{N}_c$. From the second
equality in \eqref{a1} and the definition of $J$, one obtains
\begin{equation}\label{a2}
d(c)=\dfrac{p}{2(p+1)(p+2)}c^{\frac{2}{p}-\frac{1}{2}}J(\psi),
\end{equation}
where
$\psi(x,y)=c^{-\frac{1}{p}}\ff\left(\dfrac{x}{c},\dfrac{y}{\sqrt{c}}\right)$.
Note that $\psi$ satisfies (\ref{2bozk}), with $c=1$. But, from
\eqref{2bozk-nonex-5} and \eqref{2bozk-nonex-7}, one infers that
\[
\dfrac{1}{(p+1)(p+2)}J(\ff)=\dfrac{2c}{4-p}\mathscr{F}(\ff).
\]
Thus, from \eqref{a2} follows the formula
$$
d'(c)=c^{(\frac{2}{p}-\frac{3}{2})}\mathscr{F}(\psi),
$$
whence
\[
d''(c)=\left(\dfrac{2}{p}-\dfrac{3}{2}\right)c^{(\frac{2}{p}-\frac{5}{2})}\mathscr{F}(\psi).
\]
This proves the lemma.
 \end{proof}

A study is initiated of the behavior of $d$ in a neighborhood of
the set $\mathscr{N}_c.$

\begin{lemma}
Let $c>0$. Then, there exists  a positive number $\epsilon$
and a $C^1$-map
$\mathpzc{v}:\mathpzc{B}_\epsilon(\mathscr{N}_c)\to(0,+\infty)$
defined by
\[
\mathpzc{v}(u)=d^{-1}\left(\dfrac{p}{2(p+1)(p+2)}J(u)\right)\!,
\]
such that $\mathpzc{v}(\ff)=c$ for every $\ff\in\mathscr{N}_c$, where
\[
\mathpzc{B}_\epsilon\left(\mathscr{N}_c\right)=\left\{\ff\in\mathscr{Z}\;;\;
 \inf_{\psi\in\mathscr{N}_c}\|\ff-\psi\|_{\mathscr{Z}}<\epsilon\right\}\!.
\]
\end{lemma}
 \begin{proof}
  By definition,  $\mathscr{N}_c$ is a
bounded set in $\mathscr{Z}$. Moreover,
\[
\mathscr{N}_c\subset B(0,r)\subset\mathscr{Z},
\]
where $r=(2(p+1))^{\frac{2}{p}}I_1^{\frac{p+2}{p}}$ and $B(0,r)$ is
the ball of radius $r>0$ centered at the origin in $\mathscr{Z}$.
Let $\rho>0$ be sufficiently large  that $\mathscr{N}_c\subset
B(0,\rho)\subset\mathscr{Z}$.
 Since the function $u\mapsto J(u)$ is uniformly continuous on bounded sets,
 there exists $\epsilon>0$ such that if $u,v\in B(0,\rho)$ and
$\|u-v\|_\mathscr{Z}<2\epsilon$ then $|J(u)-J(v)|<\rho$. Considering
the neighborhoods $\mathscr{I}=(d(c)-\rho,d(c)+\rho)$ and
$\mathpzc{B}_\epsilon(\mathscr{N}_c)$ of $d(c)$ and $\mathscr{N}_c$,
respectively, we have that if
$u\in\mathpzc{B}_\epsilon(\mathscr{N}_c)$ then $J(u)\in\mathscr{I}$.
Therefore $\mathpzc{v}$ is well defined on
$\mathpzc{B}_\epsilon(\mathscr{N}_c)$ and satisfies
 $\mathpzc{v}(\ff)=c$, for all $\ff\in\mathscr{N}_c$.
 \end{proof}

Here is the crucial inequality in the study of stability.

\begin{lemma}\label{stab-lem-1}
Let $c>0$ and suppose  that $d''(c)>0$. Then for all
$u\in\mathpzc{B}_\epsilon(\mathscr{N}_c)$ and any $\ff\in\mathscr{N}_c$,
\[
\mathscr{E}(u)-\mathscr{E}(\ff)+\mathpzc{v}(u)\left(\mathscr{F}(u)-\mathscr{F}
(\ff)\right)\geq\frac{1}{4}d''(c)|\mathpzc{v}(u)-c|^2\!.
\]
\end{lemma}
 \begin{proof} For $\omega>0$, let $I_\omega$ be the
functional
\[
I_\omega(\ff)=\frac{1}{2}\int_{\rr^2}\left(\omega\ff^2+\ff\mathscr{H}\ff_x+\ff_y^2\right)\;dxdy.
\]
It follows that
\[
\mathscr{E}(u)+\mathpzc{v}(u)\mathscr{F}(u)=I_{\mathpzc{v}(u)}(u)-\dfrac{1}{(p+1)(p+2)}J(u).
\]
Let $\ff_\omega$ denote any element of $\mathscr{N}_\omega$. It is
easy to see that $J(u)=J\left(\ff_{\mathpzc{v}(u)}\right)$ because
$d(\mathpzc{v}(u))=\frac{p}{2(p+1)(p+2)}J(u)$ for
$u\in\mathpzc{B}_\epsilon(\mathscr{N}_c)$ and
$d(\mathpzc{v}(u))=\frac{p}{2(p+1)(p+2)}J\left(\ff_{\mathpzc{v}(u)}\right)$.
It thus transpires that
\[
I_{\mathpzc{v}(u)}(u)\geq I_{\mathpzc{v}(u)}\left(\ff_{\mathpzc{v}(u)}\right)\!.
\]
A Taylor expansion of $d$ around the value $c$ yields
\[\begin{split}
\mathscr{E}(u)+\mathpzc{v}(u)\mathscr{F}(u)&\geq I_{\mathpzc{v}(u)}
\left(\ff_{\mathpzc{v}(u)}\right)-\dfrac{1}{(p+1)(p+2)}J\left(\ff_{\mathpzc{v}(u)}\right)\\
&=d(\mathpzc{v}(u))\geq d(c)+\mathscr{F}(\ff)(\mathpzc{v}(u)-c)+\dfrac{1}{4}d''(c)|\mathpzc{v}(u)-c|^2\\
&=\mathscr{E}(\ff)+\mathpzc{v}(u)\mathscr{F}(\ff)+\dfrac{1}{4}d''(c)|\mathpzc{v}(u)-c|^2\!,
\end{split}\]
and the lemma follows.
 \end{proof}

Before proving stability, we state a well-posedness
result for (\ref{2bozk-main}).  This can be proved in several standard ways,
for example by using a
parabolic regularization  (see \cite{iorio} and \cite{CP}).
\begin{theorem}
 Let $s>2$. Then for any $u_0\in H^s(\rr^2)$, there exist
$T=T(\|u_0\|_{H^s})>0$ and  a unique solution $u\in
C([0,T];H^s(\rr^2))$ of   equation (\ref{2bozk-main}) with
$u(0)=u_0$. In addition, $u(t)$ depends continuously on $u_0$  in
the $H^s-$norm and satisfies $\mathscr{E}(u(t))=\mathscr{E}(u_0)$,
$\mathscr{F}(u(t))=\mathscr{F}(u_0)$, for all $t\in[0,T)$.
\end{theorem}

When $0 < p < \frac43$, the stability  in $\mathscr{Z}$ of the set of minimizers
$\mathscr{N}_c$ is established next.

\begin{theorem}
Let $c>0, \, s > 2, \, 0<p< \frac43$ and $\lam=\big(2(p+1)I_1\big)^\frac{p+2}{p}$. Then the
set $\mathscr{N}_c=M_\lam$ is $\mathscr{Z}$-stable with regard to
the flow of the BO-ZK equation. That is, for any positive $\epsilon$, there is a positive
$\delta = \delta(\epsilon)$ such that if $u_0\in H^s$ and
$\inf_{\ff \in \mathscr{N}_c} \|u_0-\ff\|_{H^s}\leq\delta$, then the solution $u(t)$ of
(\ref{2bozk-main}) with $u(0)=u_0$ satisfies
\[
\sup_{t\geq0}\inf_{\psi\in\mathscr{N}_c}\|u(t)-\psi\|_\mathscr{Z}\leq\epsilon.
\]
\end{theorem}
 \begin{proof} Assume that $\mathscr{N}_c$ is
$\mathscr{Z}$-unstable with regard to the flow of the BO-ZK
equation. Then, there is a sequence of initial data
$u_k(0)\in
H^s\left(\rr^2\right)$ such that
\begin{equation}  \label{temp12}
\inf_{\ff\in\mathscr{N}_c}\|u_k(0)-
\ff\|_{H^s} \leq \frac{1}{k} \quad\; {\rm and} \quad
\sup_{t\in[0,T)}\inf_{\psi\in\mathscr{N}_c}\|u_k(t)-
\psi\|_\mathscr{Z}\geq\epsilon,
\end{equation}
where $u_k(t)$ is the solution of (\ref{2bozk-main}) with initial
data $u_k(0)$. By continuity in $t$, for all $k$
large enough, there are times $t_k$ such that
\begin{equation}
\inf_{\ff\in\mathscr{N}_c}\|u_k\left(t_k\right)-
\ff\|_\mathscr{Z}=\frac{\epsilon}{2}.\label{contr-stab}
\end{equation}
 Since $\mathscr{E}$ and $\mathscr{F}$ are
conserved quantities, it follows from \eqref{temp12}
that
\begin{gather} \label{zero}
|\mathscr{E}(u_k(t_k))-\mathscr{E}(\ff_k)|=|\mathscr{E}(u_k(0))-\mathscr{E}(\ff_k)|\to0,\\
|\mathscr{F}(u_k(t_k))-\mathscr{F}(\ff_k)|=|\mathscr{F}(u_k(0))-
\mathscr{F}(\ff_k)|\to0, \label{one}
\end{gather}
as $k\to+\infty$.  In this circumstance, Lemma \ref{stab-lem-1} implies that
\[
\mathscr{E}(u_k(t_k))-\mathscr{E}(\ff_k)+\mathpzc{v}(u_k(t_k))\big(\mathscr{F}(u_k(t_k))-
\mathscr{F}(\ff_k)\big)
\geq\frac{1}{4}d''(c)|\mathpzc{v}(u_k(t_k))-c|^2\!,
\]
for all $k$ large enough. Since $\{u_k(t_k)\}$ is
uniformly bounded in $k$, the right-hand side of the last inequality goes to zero as
$k \to +\infty$ on account of \eqref{zero} and \eqref{one}.  This in turn implies that
$\mathpzc{v}(u_k(t_k))\to c$ as $k\to+\infty$. Hence, by the definition of $\mathpzc{v}$
and continuity of $d$, we must have
\begin{equation}
\lim_{k\to+\infty}J(u_k(t_k))=\dfrac{2(p+1)(p+2)}{p}\;d(c).\label{equal-lim-1}
\end{equation}
On the other hand,  Lemma \ref{stabl-lem-2} implies that
\[\begin{split}
I(u_k(t_k))&=\mathscr{E}(u_k(t_k))+c\mathscr{F}(u_k(t_k))+\dfrac{1}{(p+1)(p+2)}J(u_k(t_k))\\
&=d(c)+ \mathscr{E}(u_k(t_k))-E(\ff_k)+c\left(\mathscr{F}(u_k(t_k))-
\mathscr{F}(\ff_k)\right) \\ \qquad \qquad \qquad & \qquad\quad+\dfrac{1}{(p+1)(p+2)}J(u_k(t_k)).
\end{split}\]
The limit (\ref{equal-lim-1}) then yields
\begin{equation}
\lim_{k\to+\infty}I(u_k(t_k))=\dfrac{p+2}{p}d(c)= \big(2(p+1)\big)^\frac{2}{p}\;I_1^{\frac{p+2}{p}}.
\label{equal-lim-2}
\end{equation}
Defining
\[
\vartheta_k(t_k)=\Big(J\big(u_k(t_k)\big)\Big)^{-\frac{1}{p+2}}u_k(t_k),
\]
it is seen that $J\left(\vartheta_k(t_k)\right)=1$.  Combining
 (\ref{equal-lim-1}), (\ref{equal-lim-2}) and Lemma
\ref{stabl-lem-2} leads to
\begin{equation}
\lim_{k\to+\infty}I(\vartheta_k(t_k))=I_1.
\end{equation}
Hence $\left\{\vartheta_k(t_k)\right\}$ is a minimizing sequence for
$I_1$. Thus, from Theorem  \ref{stab-theo1}, there exists a sequence
$\{\psi_k\}\subset M_1$ such that
\begin{equation}
\lim_{k\to+\infty}\left\|\vartheta_k(t_k)-\psi_k\right\|_\mathscr{Z}=0.\label{equal-lim-3}
\end{equation}
The Lagrange multiplier theorem then implies there is a sequence
$\{\theta_k\} \subset\rr$ such that
\begin{equation}
\mathscr{H}\left(\psi_k\right)_x+c\psi_k-\left(\psi_k\right)_{yy}=\theta_k(p+2)\psi_k^{p+1}.
\end{equation}
In other words, $2I_1=\theta_k(p+2)$, which implies
$\theta_k=\theta$ for all $k$.  Write $\ff_k=\mu\psi_k$ with
$$\mu^p=\theta(p+1)(p+2)=2(p+1)I_1.  \vspace{.11cm}
$$
Then the $\ff_k$ satisfy
(\ref{2bozk}) and $2(p+1)I(\ff_k)=J(\ff_k)=\mu^{p+2}$ so
that $\ff_k\in\mathscr{N}_c$ for all $k$.   Additionally,
(\ref{equal-lim-1})-(\ref{equal-lim-3}) and Lemma \ref{stabl-lem-2}
together allow the conclusion
\[\begin{split}
\|u_k(t_k)&-\ff_k\|_{H^s}=J\big (u_k(t_k)\big)^\frac{1}{p+2}
\left\|J\big(u_k(t_k)\big)^{-\frac{1}{p+2}}\big(u_k(t_k)-\ff_k\big)\right\|_{H^s}\\
&\leq J\big(u_k(t_k)\big)^\frac{1}{p+2}\left(\left\|\vartheta_k(t_k)-\mu^{-1}\ff_k
\right\|_{H^s}+\mu^{-1}\|\ff_k\|_{H^s}-J\big(u_k(t_k)\big)^{-\frac{1}{p+2}}\right)\!.
\end{split}\]
This in turn implies that \[
\lim_{k\to+\infty}\|u_k(t_k)-\ff_k\|_\mathscr{Z}=0,
\]
which  contradicts (\ref{contr-stab}) and  completes the proof of
the Theorem.
 \end{proof}

\section{Decay and Regularity}\label{decay and regul}

To investigate the regularity and the spatial asymptotics
of the solitary-wave solutions of (\ref{2bozk-main}), it is
convenient to take the Fourier transform of equation \eqref{2bozk}
for the solitary-wave in both $x$
and $y$.  If $(\xi,\eta)$ are the variables dual to $(x,y)$ by way of
the Fourier transform, then
\eqref{2bozk} implies that
\begin{equation} \label{transformed}
\widehat{\ff}=\dfrac{\widehat{\mathpzc{g}}}{c-\alpha|\xi|+\varepsilon\eta^2}, \quad {\rm where}
\quad
\mathpzc{g}\, =\, -\frac{1}{p+1}\ff^{p+1}\!.
\end{equation}
Taking the inverse Fourier transform then yields
\begin{equation} \label{transformedeqn}
\ff \, = \, -\frac{1}{p+1}\int_{\rr^2 } K\big(x- s,y - t\big) \ff^{p+1}(s,t)\,dsdt.
\end{equation}

Properties of
the  integral kernel $K$ in  (\ref{transformedeqn}) will be central in the analysis to follow.
   Here are few standard
properties of anisotropic Sobolev spaces that will be helpful in expressing
useful aspects of  $K$.
\begin{lemma}\label{algebra}
If $s_i>1/2$, for  $i=1,2$, then
 $H^{s_1,s_2}$ is an algebra.
\end{lemma}

\begin{lemma}\label{interp-aniso} Let $s_{ij}, 1 \leq i,j \leq 2$ and $\theta \in [0,1]$
be given real numbers with $s_{1j}\leq  s_{2j}$, $j=1,2$.  Define $\varrho_j = \theta s_{1j} + (1-\theta)s_{2j}$
for $j=1,2$.  Then, $H^{\varrho_1,\varrho_2}$ is an interpolation space between
$H^{s_{11},s_{12}}$ and it's subspace $H^{s_{21},s_{22}}$.  Moreover, if $f \in H^{s_{21},s_{22}}$,
then
\begin{equation}
\|f\|_{H^{\varrho_1,\varrho_2}}\leq
\|f\|_{H^{s_{11},s_{12}}}^\theta\;\|f\|_{H^{s_{21},s_{22}}}^{1-\theta}.
\end{equation}
\end{lemma}

\begin{remark}
Since $\widehat{K}(\xi,\eta)=\dfrac{1}{c-\alpha|\xi|+\eta^2}$,  the Residue Theorem
allows us to write the kernel $K$ as an integral, namely
\begin{equation}
K(x,y)= K_c(x,y)= C\int_0^{+\infty}\dfrac{|\alpha|\sqrt{t}}{\alpha^2t^2+x^2}\;
e^{-\left(ct+\frac{y^2}{4t}\right)}\;dt,\label{kernel}
\end{equation}
where $C>0$ is  independent of $\alpha$, $x$ and $y$.
 Fubini's theorem can then be used to show that
\[
\|K\|_{L^1}=C\int_0^{+\infty}\int_{\rr^2}\dfrac{|\alpha|\sqrt{t}}{\alpha^2t^2+x^2}\;
e^{-\left(ct+\frac{y^2}{4t}\right)}dxdydt=C(\alpha)\int_0^{+\infty}e^{-ct}dt.
\]
\end{remark}
\noindent  In consequence of   representation \eqref{kernel}, the following facts
about $K$ become clear.

\begin{lemma}\label{prop1}  The kernel
$K$ is positive, an even function of both $x$ and $y$, monotone decreasing
in both  $|x|$ and $|y|$, tends to zero as $|(x,y)| \to \infty$  and is $C^\infty$
away from the origin.   Furthermore,
 $\widehat{K}\in L^p(\rr^2)$ for any $p\in(3/2,+\infty]$ and
 $K\in L^p(\rr^2)$, for any $p\in[1,3)$.  (However, while $K(x,y)$
is symmetric in both $x$ and $y$, it is not radially symmetric.)
\end{lemma}

\begin{lemma}\label{prop2}
$K\in H^{s_1,0}\left(\rr^2\right)\cap H^{0,s_2}\left(\rr^2\right)$
for any $s_1<\frac{1}{4}$ and $s_2<\frac{1}{2}$. Moreover,
$K\in H^{r,s}\left(\rr^2\right)\cap H^{s_1,s_2}\left(\rr^2\right)$,
where $rs_2+ss_1=s_1s_2$ and $r\in[0,1]$.
\end{lemma}
\begin{lemma} \label{prop3}
\begin{enumerate}[(i)]
\item  $\widehat{K}\in H^{s_1,0}\left(\rr^2\right)\cap H^{0,s_2}
\left(\rr^2\right)$, for any $s_1<\frac{3}{2}$ and $s_2\in\rr$.
Moreover, $\widehat{K}\in H^{r,s}\left(\rr^2\right)\cap H^{(s_1,s_2)}\left(\rr^2\right)$, where $rs_2+ss_1=s_1s_2$ and $r\in[0,1]$.
\item  $\widehat{K}\in H^{(s_1,s_2)}_p\left(\rr^2\right)$,
 for any $s_1<1+\frac{1}{p}$, $p\geq2$ and $s_2\in\rr$.
 \item   $|x|^{s_1}|y|^{s_2}K\in L^p\left(\rr^2\right)$, for any $s_1,s_2\geq0$ and $1\leq p\leq\infty$ such that
  $s_1<2-\frac{1}{p}$ and $2s_1+s_2>1-\frac{3}{p}$.
\end{enumerate}
\end{lemma}

With these facts about $K$ in hand,  the solitary-wave solutions of
the BO-ZK equation \eqref{2bozk-main} now become the focus of attention.

\begin{theorem}\label{regularity}  Let $p$ be a positive integer.
Any solitary-wave solution $\ff$ of (\ref{2bozk-main})
belongs to $H^{(k)}_r$, for all $k\in\N$ and all
$r\in[1,+\infty]$. In particular, the solitary-wave solutions of the
BO-ZK equation are
continuous, bounded and tend to zero at infinity.
\end{theorem}
 \begin{proof}
  Formula \eqref{transformed}
implies that $\ff\in H^{\frac{1}{2},1}(\rr^2)\cap
H^{0,2}(\rr^2)\cap H^{1,0}(\rr^2)$.  Lemma
\ref{interp-aniso} and the embedding (\ref{embedding-0}) then imply
that $\ff\in H^{s,2(1-s)}(\rr^2)$, for any $s\in[0,1]$.  A
bootstrapping argument and the use of  Lemmas \ref{interp-aniso} and
 \ref{algebra} completes the proof.
 \end{proof}

More detailed aspects of the solitary-wave solutions of \eqref{2bozk-main}
 are now addressed.  Interest will focus first upon their  symmetry properties.
   For $u:\rr^2\to\rr^+$, $u^\sharp$ will
denote the Steiner symmetrization of $u$ with respect to $\{x=0\}$ and $u^*$ the Steiner symmetrization of $u$ with respect to $\{y=0\}$ (see, for example, \cite{bll,k1,vans}).
Notice that $u^{\sharp*} = u^{*\sharp}$ is a function symmetric with respect to both
the $x$- and $y$-axis.

\begin{lemma} \label{prop4}
If  $f\in\mathscr{Z}$, then $ |f|$ lies in $\mathscr{Z}$ and $I(|f|) \leq I(f)$.
\end{lemma}
 \begin{proof}
   If $g=|f|$, then for any $c > 0$,
\[
\left\langle f,K\ast f\right\rangle\leq\left\langle g,K\ast g\right\rangle\!,
\]
where $K = K_c$.  It thus transpires that
\[
\int_{\rr^2}\what{K}(\xi,\eta)\left|\what{f}(\xi,\eta)\right|^2\!d\xi d\eta=\left\langle f,K\ast f\right\rangle
\leq\left\langle g,K\ast g\right\rangle=\int_{\rr^2}\what{K}(\xi,\eta)\left|\what{g}(\xi,\eta)\right|^2\!d\xi d\eta.
\]
Since  $\big\|\what{f}\big\|_{L^2}=\big\|\what{g}\big\|_{L^2}$, it follows that
\begin{equation}
\int_{\rr^2}c\left(1-c\what{K}\right)\;|\what{g}(\xi,\eta)|^2\;
d\xi d\eta\leq\int_{\rr^2}c\left(1-c\what{K}\right)\;\left|\what{f}
(\xi,\eta)\right|^2\!d\xi d\eta. \label{inq-pos-kernel}
\end{equation}
Taking the limit as $c\to+\infty$ on both sides of
(\ref{inq-pos-kernel}), the Monotone Convergence Theorem yields
\begin{equation}   \label{decreasing}
\int_{\rr^2}\left(|\xi|+\eta^2\right)|\widehat{g}(\xi,\eta)|^2
d\xi d\eta\leq\int_{\rr^2}\left(|\xi|+\eta^2\right)\left|\what{f}(\xi,\eta)\right|^2\!d\xi d\eta,
\end{equation}
which shows that $|f|\in\mathscr{Z}$ and that $I(|f|) \leq I(f)$.
 \end{proof}

\begin{corollary}
For $c> 0$, there is always a non-negative solitary-wave solution $\ff_c$
of the BO-ZK equation.
\end{corollary}
\begin{proof}  Theorem \ref{varitional-charac} assures that there are solitary-wave solutions $\psi$, say.
The last result shows that if $\psi \in M_\lam$, then so is $\ff = |\psi|$.
\end{proof}

If $p = \frac{k}{m}$ where $m$ is odd and $k$ and $m$ relatively prime it follows from the formula
\begin{equation}\label{conv-form}
\ff = \frac{1}{p+1} K*\ff^{p+1}
\end{equation}
that if $k$ is odd, then necessarily all solitary-wave solutions are non-negative.  This is false if
$k$ is even, however.  Indeed, in this case, if $\ff$ is a solitary wave, then so is $-\ff$.   Hence,
when $k$ is even, there are always at least two solitary-wave solutions, one
positive and one negative.  Of course, when $k$ is even, it is also the case that
$J(|f|) = J(f)$.

\begin{lemma} \label{sym} If $f \in \mathscr{Z}$ is non-negative, it's Steiner symmetrizations  $f^\sharp$ and $f^*$
also lie in $\mathscr{Z}$. Moreover, $I(f^\sharp) \leq I(f)$ and $I(f^*) \leq I(f)$.
\end{lemma}
\begin{proof}
 Remark first that $K^\sharp=K=K^*$. The Reisz-Sobolev rearrangement inequality (see \cite{bll,k1,vans}) implies that
\[\begin{split}
\int_{\rr^4}&f(x,y)f(s,t)K(x-s,y-t)ds\,dt\,dx\,dy\\
&\qquad\qquad\leq\int_{\rr^4}f^\sharp(x,y)f^\sharp(s,t)K(x-s,y-t) ds\,dt\,dx\,dy.
\end{split}\]
In the Fourier transformed variables, this amounts to
\[
\int_{\rr^2}\what{K}(\xi,\eta)\left|\what{f}(\xi,\eta)\right|^2\! d\xi d\eta \leq
\int_{\rr^2}\what{K}(\xi,\eta)\left|\what{f^\sharp}(\xi,\eta)\right|^2\! d\xi d\eta.
\]
On the other hand, the fact that symmetrization does not change the
measure theoretic properties of $f$ implies that
\[
\big\|\what{f}\big\|_{L^2(\rr^2)}=\big\|f\big\|_{L^2(\rr^2)}
=\big\|f^\sharp\big\|_{L^2(\rr^2)}=\big\|\what{f^\sharp}\big\|_{L^2(\rr^2)}.
\]
 This together with the analysis in Lemma \ref{prop4} shows that $f^\sharp\in\mathscr{Z}$
 and that  $I(f^\sharp) \leq I(f)$.
The same argument applies to
$f^*$.
\end{proof}

Since Steiner symmetrization preserves the $L^{p+2}-$norm, it follows that $J(\ff)=J(\ff^\sharp)$. In consequence of Lemma \ref{sym},
\[
I_\lam\leq I\left(\ff^\sharp\right)\leq I(\ff)=I_\lam.
\]
Therefore $\ff^\sharp\in M_\lam$.  The same argument shows that  $\ff^*\in M_\lam$.
\hfill$\Box$

\begin{corollary}\label{symmetry}
There are non-negative, solitary-wave solutions of the BO-ZK equation \eqref{2bozk-main} that are  symmetric
with respect to both the propagation direction and the transverse direction and are  monotone decreasing in both $|x|$ and $|y|$.
\end{corollary}
\begin{proof}
By Theorems \ref{existence} and \ref{regularity}, there is a
non-negative function $\ff $ satisfying \eqref{2bozk}.
Since Steiner symmetrization preserves the $L^{p+2}-$norm, it follows that
$J(\ff)=J(\ff^\sharp) = J(\ff^{\sharp*})$. On the other hand, because
of Lemma \ref{sym},
the double rearrangement $\ff^{\sharp*}$   has the property that
\[
I_\lam\leq I\left(\ff^{\sharp*}\right)\leq I(\ff^\sharp) \leq I(\ff)=I_\lam.
\]
Therefore, $\ff^{\sharp*}$
 is a non-negative solitary-wave solution of  equation \eqref{2bozk-main} which is symmetric with respect to both $\{x=0\}$ and $\{y=0\}$ and which is monotone decreasing with respect to both $|x|$
and $|y|$.
 \end{proof}

\begin{remark}
One may also obtain symmetry properties of the solitary-wave solutions of \eqref{2bozk-main} by using the reflection method and a unique continuation argument (see \cite{lopes-maris} and \cite{prep-3}).
\end{remark}

\section{Spatial Asymptotics}\label{decaysec}

Attention is now turned to the spatial decay properties of the solitary-wave solutions of
 \eqref{2bozk-main}.   In this analysis, we follow the lead of \cite{bonali}.

\begin{lemma}\label{decay-integ}
Let $j\in\mathbb{N}$. Suppose also that $\ell$ and $m$ are two constants satisfying $0<\ell<m-j$.
Then there exists $C>0$, depending only on $\ell$ and $m$,
such that for all $\epsilon\in(0,1]$,  we have
\begin{equation}\label{est-1}
\int_{\rr^j}\frac{|a|^\ell}{(1+\epsilon|a|)^m(1+|b-a|)^m}\;\dd a
\leq\frac{C\;|b|^\ell}{(1+\epsilon|b|)^m},\qquad\forall\;b\in\rr^j,\;|b|\geq1,
\end{equation}
and
\begin{equation}\label{est-2}  \hspace{-1.3cm}
\int_{\rr^j}\frac{\dd a}{(1+\epsilon|a|)^m(1+|b-a|)^m}
\leq\frac{C}{(1+\epsilon|b|)^m},\qquad\forall\;b\in\rr^j.
\end{equation}
\end{lemma}
The proof of this elementary lemma is
 essentially the same as the proof of Lemma 3.1.1 in \cite{bonali} (see \cite{jpha}).

\begin{theorem}\label{decay-theo-1}
Let $\ff$ be a solitary-wave solution of (\ref{2bozk}).
\begin{enumerate}[(i)]
\item  For all $q\in(3/2,+\infty)$, $\ell\in[0,1)$  $\varrho\geq0$,
$|x|^\ell|y|^\varrho\ff(x,y)\in L^q\left(\rr^2\right)$.

\item   For all $q\in(3/2,+\infty)$ and any $\theta\in[0,1)$, $|(x,y)|^\theta\ff(x,y)\in L^q\left(\rr^2\right)$.

\item  And finally, $\ff\in L^1\left(\rr^2\right)$.
\end{enumerate}
\end{theorem}
 \begin{proof}
  \rm{(i)}\quad For $q\in(1,3)$ and  $1-\frac{1}{q}<s_1<2-\frac{1}{q}$,
let $\ell\in\left[0,s_1-1+\frac{1}{q}\right)$. Also, for $s_2>1-\frac{1}{q}$, choose $\varrho\in\left[0,s_2-1+\frac{1}{q}\right)$. For $0<\epsilon<1$, define
$\mathpzc{h}_\epsilon$ by
\[
\mathpzc{h}_\epsilon(x,y)=\mathpzc{A}(x,y)\;\ff(x,y),
\]
where
\[
\mathpzc{A}(x,y)=\frac{|x|^\ell|y|^\varrho}{(1+\epsilon|x|)^{s_1}(1+\epsilon|y|)^{s_2}},
\]
 Then, by using the explicit representation of $\mathpzc{h}_\epsilon$, it is easy to check that
$\mathpzc{h}_\epsilon\in L^{q'}\left(\rr^2\right)$, where $q'=\frac{q}{q-1}$.
H\"{o}lder's inequality and \eqref{transformedeqn} then implies that
\[
|\ff(x,y)|\leq C(s_1,s_2,q)\left(\int_{\rr^2}|\mathpzc{G}_{x,y}(z,w)|^{q'}\;dzdw\right)^\frac{1}{q'},
\]
where
\[
\mathpzc{G}_{x,y}(z,w)=\frac{\mathpzc{g}(\ff)(z,w)}
{\big(1+|x-z|\big)^{s_1}\big(1+|y-w|\big)^{s_2}},
\]
$\mathpzc{g}(t)=\frac{t^{p+1}}{p+1}$ and
 $$
 C:=C(s_1,s_2,p)=\big\|(1+|x|)^{s_1}(1+|y|)^{s_2}
K\big\|_{L^q\left(\rr^2\right)}<\infty.
$$
This last constant is finite thanks to Lemma \ref{prop3}. Since the solitary wave $\ff$ converges to the rest state as $|(x,y)|\to+\infty$, it follows that
for every $\delta>0$, there exists $R_\delta>1$
 such that if $|(x,y)|\geq R_\delta$, then
\[
\left|\mathpzc{g}(\ff)(x,y)\right|\leq\delta|\ff(x,y)|.
\]

An application of H\"{o}lder's inequality yields
\[\begin{split}
\int_{\rr^2\setminus B(0,R_\delta)}|\mathpzc{h}_\epsilon&(x,y)|^{q'}\!
dxdy=\int_{\rr^2\setminus B(0,R_\delta)}|\mathpzc{h}_\epsilon(x,y)|^{q'-r}
\mathpzc{A}^r(x,y)|\ff(x,y)|^rdxdy\\
&\leq C^r\int_{\rr^2\setminus B(0,R_\delta)}|\mathpzc{h}_\epsilon(x,y)|^{q'-r}
\mathpzc{A}^r(x,y)\left\|\mathpzc{G}_{x,y}\right\|_{L^{q'}(\rr^2)}^r(x,y)\,dxdy\\
&\leq C^r\|\mathpzc{h}_\epsilon\|_{L^{q'}\left(\rr^2\setminus
B(0,R_\delta)\right)}^{q'-r}\left\|\;\mathpzc{A}\;\|\mathpzc{G}_{x,y}\|_{L^{q'}
(\rr^2)}\right\|_{L^{q'}\left(\rr^2\setminus B(0,R_\delta)\right)}^r.
\end{split}\]
Because $\mathpzc{h}_\epsilon\in L^{q'}\left(\rr^2\right)$, the latter inequality implies
\[
\|\mathpzc{h}_\epsilon\|_{L^{q'}\left(\rr^2\setminus
B(0,R_\delta)\right)}^{r}\leq C^r \left\|\;\mathpzc{A}\;\|\mathpzc{G}_{x,y}\|_{L^{q'}
(\rr^2)}\right\|_{L^{q'}\left(\rr^2\setminus B(0,R_\delta)\right)}^r,
\]
 which is to say,
\[
\int_{\rr^2\setminus B(0,R_\delta)}|\mathpzc{h}_\epsilon(x,y)|^{q'}\!dxdy
\leq C^{q'}\int_{\rr^2\setminus B(0,R_\delta)}\mathpzc{A}^{q'}(x,y)\|
\mathpzc{G}_{x,y}\|^{q'}_{L^{q'}(\rr^2)}dxdy.
\]
Fubini's theorem and Lemma \ref{decay-integ} combine to show that
\begin{equation}\label{mideeq}
\begin{split}
&\int_{\rr^2\setminus B(0,R_\delta)}\mathpzc{A}^{q'}(x,y)\|
\mathpzc{G}_{x,y}\|^{q'}_{L^{q'}(\rr^2)}(x,y)\,dxdy\\
&=\int_{\rr^2}\left|\mathpzc{g}(\ff)(z,w)\right|^{q'}\left(\int_{\rr^2\setminus
B(0,R_\delta)}\frac{\mathpzc{A}^{q'}(x,y)}{(1+|x-z|)^{q's_1}(1+|y-w|)^{q's_2}}dxdy\right)\!dzdw\\
&\leq {C} \int_{\rr^2\setminus B(0,R_\delta)}\left|\mathpzc{g}
(\ff)(z,w)\right|^{q'}\mathpzc{A}^{q'}(z,w)\,dzdw\\
&+ \int_{B(0,R_\delta)}\left|\mathpzc{g}(\ff)(z,w)\right|^{q'}
\left(\int_{\rr^2\setminus B(0,R_\delta)}\frac{\mathpzc{A}^{q'}(x,y)}{(1+|x-z|)^{q's_1}(1+|y-w|)^{q's_2}}dxdy\right)\!dzdw,
\end{split}
\end{equation}
where we used \eqref{est-1} (with $j=1$) to show that for $|(z,w)|$ large,
\[
\int_{\rr^2\setminus
B(0,R_\delta)}\frac{\mathpzc{A}^{q'}(x,y)}{(1+|x-z|)^{q's_1}(1+|y-w|)^{q's_2}}dxdy\leq C \mathpzc{A}^{q'}(z,w).
\]
The second integral on the right-hand side of \eqref{mideeq} is bounded by a constant, say ${C}'$, depending
 on $\ff$ and $R_\delta$, but independent of $\epsilon$. Therefore,
 by using the fact that $\left|\mathpzc{g}(\ff)(x,y)\right|\leq\delta|\ff(x,y)|$
 on $\rr^2\setminus B(0,R_\delta)$, there obtains
\[
\int_{\rr^2\setminus B(0,R_\delta)}|\mathpzc{h}_\epsilon(x,y)|^{q'}\!
dxdy\leq C^{q'}\left({C}\delta^{q'}\int_{\rr^2\setminus
B(0,R_\delta)}|\mathpzc{h}_\epsilon(x,y)|^{q'}\;\!dxdy+{C}'\right).
\]
Choosing $\delta$ such that $C\delta{C}^{\frac{1}{q'}}<1$,
 the last inequality entails that
\begin{equation}
\int_{\rr^2\setminus B(0,R_\delta)}|\mathpzc{h}_\epsilon(x,y)|^{q'}\!
dxdy\leq{C^{''}}\!,\label{decay--ineq-1}
\end{equation}
where ${C^{''}}$ is a constant independent of $\epsilon$.
Letting $\epsilon\to0$ in (\ref{decay--ineq-1}) and applying  Lebesgue's
dominated convergence theorem,  one deduces
\[
\int_{\rr^2\setminus B(0,R_\delta)}|x|^{\ell q'}|y|^{\varrho q'}|\ff(x,y)|^{q'}\!dxdy\leq{C}.
\]
Hence $|x|^\ell|y|^\varrho\ff(x,y)\in L^{q'}\left(\rr^2\right)$,
for $q'=\frac{q}{q-1}$.\\\indent
In  the limits $q\to1$ and $q\to3$, we have
 $\ell\to1$ and $q'\in(3/2,+\infty)$. This proves part (i) of the theorem.

\rm{(ii)}\quad This follows directly from \rm{(i)}.

\indent \rm{(iii)}\quad Let $s>1$ and $\mathpzc{g}$, $\delta$
and $R_\delta$  be as defined in the proof of   \rm{(i)}. For $\epsilon > 0$
let $\mathpzc{A}$
be
\[
\mathpzc{A}_\epsilon(x,y)=\frac{1}{\left(1+\epsilon|(x,y)|\right)^{s}}.
\]
 Fubini's Theorem,  Lemma \ref{decay-integ} and the fact
that $\ff$, \!$\mathpzc{A}_\epsilon\in L^2\left(\rr^2\right)$ so that the product
$\ff\mathpzc{A}_\epsilon\in L^1\left(\rr^2\right)$ allow us
to adduce the inequalities
\[\begin{split}
\int_{\rr^2\setminus B\left(0,R_\delta\right)}&|\ff(x,y)|\mathpzc{A}_\epsilon(x,y)\,dxdy\\
&\leq\int_{\rr^2}\left|\mathpzc{g}(\ff)(z,w)\right|\left(\int_{\rr^2\setminus B\left(0,R_\delta\right)}\mathpzc{A}_\epsilon(x,y)K(x-z,y-w)\,dxdy\right)\!dzdw\\
&\leq\int_{\rr^2}\left|\mathpzc{g}(\ff)(z,w)\right|\left(\int_{\rr^2\setminus B\left(0,R_\delta\right)}\mathpzc{A}_1^{-2}(x-z,y-w)K^2(x-z,y-w)\,dxdy\right)^\frac{1}{2}\\
&\qquad\times\left(\int_{\rr^2\setminus B\left(0,R_\delta\right)}\mathpzc{A}_1^2(x-z,y-w)\mathpzc{A}_\epsilon^2(x,y)\;dxdy\right)^\frac{1}{2}\!dzdw\\
&\leq C(s){C}^\frac{1}{2}\int_{\rr^2}\left|\mathpzc{g}(\ff)(z,w)\right|\mathpzc{A}_\epsilon(z,w)\,dzdw\\
&\leq C(s){C}^\frac{1}{2}\delta\int_{\rr^2\setminus B\left(0,R_\delta\right)}|\ff(z,w)|\mathpzc{A}_\epsilon(z,w)\,dzdw\\
&\qquad+C(s){C}^\frac{1}{2}\int_{B(0,R_\delta)}
\left|\mathpzc{g}(\ff)(z,w)\right|dzdw.
\end{split}\]
Letting $\epsilon \to 0$, Fatou's lemma together with the restriction on $\delta$ leads to the conclusion that $\ff\in L^1\left(\rr^2\right)$.
\end{proof}

 Theorem \ref{decay-theo-1}, identity
 \eqref{conv-form} and  the elementary
inequality
\begin{equation}
|t|^\theta\leq C\left(|t-s|^\theta+|s|^\theta\right),\quad\mbox{for}\;\;\theta\geq0.\label{conv-classic}
\end{equation}
imply the following.

\begin{corollary}\label{infinity-decay}
Suppose that $\ff\in L^\infty\left(\rr^2\right)$ satisfies
(\ref{2bozk}) and $\ff\to0$ at infinity. Then
\begin{enumerate}[(i)]
\item  $|x|^\ell|y|^\varrho\ff(x,y)\in L^\infty\left(\rr^2\right)$,
for all  $\ell\in[0,1)$ and any $\varrho\geq0$,
\item  $|(x,y)|^\theta\ff(x,y)\in L^\infty\left(\rr^2\right)$, for all  $\theta\in[0,1)$.
\end{enumerate}
\end{corollary}

The aim now is to display even stronger decay properties in the $x$-variable
 for solitary-wave solutions of  the {\it BO-ZK} equation.
These results are developed in a sequence of lemmas.

\begin{lemma}\label{lem-infin-max}
$|x|^2|y|^\varrho K\in L^\infty\left(\rr^2\right)$, for any
  $\varrho\geq0$.
\end{lemma}
\begin{proof}
In view of the explicit form of $K$, the proof is straightforward.
\end{proof}

\begin{corollary}
$|x|^\ell|y|^\varrho\ff(x,y)\in L^\infty\left(\rr^2\right)$, for any
$0\leq\ell\leq2$ and any $\varrho\geq0$.
\end{corollary}
 \begin{proof}
The proof is based on a standard bootstrapping argument.  Decay in the
$y$-variable is not in question, so without loss of generality, take it that
that $\varrho=0$. Setting $\gamma_1=\min\{2,p+1\}$ and making use of the inequality
\begin{equation}
|x|^{\gamma_1}|\ff|\lesssim|x| ^{\gamma_1}|K|\ast|\mathpzc{g}(\ff)|+|K|\ast||x|^{\gamma_1}|\mathpzc{g}(\ff)||,
\end{equation}
where $\mathpzc{g}(t)=\frac{t^{p+1}}{p+1}$,
we obtain from Corollary \ref{infinity-decay},
Lemma \ref{lem-infin-max} and Theorem \ref{decay-theo-1} that $|x|^{\gamma_1}\ff\in L^\infty(\rr^2)$.   The proof is compete if $\gamma_1=2$. If $\gamma_1<2$, then  define $\gamma_2=\min\{2,(p+1)^2\}$ and repeat the above argument to show
$|x|^{\gamma_2}\ff\in L^\infty(\rr^2)$.  Continuing in this manner, one concludes that
 $|x|^2\ff\in L^{\infty}(\rr^2)$ after a finite number of steps.
 \end{proof}

The following corollary follows from (\ref{conv-classic}),
Corollary \ref{infinity-decay} and Theorem \ref{decay-theo-1}.
\begin{corollary}
\begin{enumerate}[(i)]
\item  $|x|^\ell|y|^\varrho\ff(x,y)\in L^1\left(\rr^2\right)$, for all  $\ell\in[0,1)$ and any $\varrho\geq0$,
\item  $|(x,y)|^\theta\ff(x,y)\in L^1\left(\rr^2\right)$, for all  $\theta\in[0,1)$.
\end{enumerate}
\end{corollary}
\begin{lemma}\label{decay-expon-lem}
For any $1\leq r,q<\infty$, there is $\sigma_0>0$ such that for all $\sigma\in[0,\sigma_0)$ and $s\in(\frac{1}{2}-\frac{1}{r}-\frac{1}{2q},2-\frac{1}{r})$, we have
\[
|x|^se^{\sigma|y|}K\in L^r_xL_y^q(\rr^2)\cap L_y^qL^r_x(\rr^2).
\]
\end{lemma}
\begin{proof}
It suffices to choose $\sigma_0=\sqrt{\frac{c}{q}}$, where $c$ is the wave velocity and use \eqref{kernel}.
\end{proof}

The next result  is a consequence of another of Young's inequalities, namely
\[
\|f\ast g\|_{L_y^qL_x^r(\rr^2)}\leq\|f\|_{L_y^{q_1}L_x^{r_1}(\rr^2)}\|g\|_{L_y^{q_2}L_x^{r_2}(\rr^2)},
\]
where $1\leq r,q,r_1,q_1,r_2,q_2\leq\infty$,
$1+\frac{1}{r}=\frac{1}{r_1}+\frac{1}{r_2}$ and
 $1+\frac{1}{q}=\frac{1}{q_1}+\frac{1}{q_2}$.

\begin{corollary}   \label{correg}
$\ff\in L^r_xL_y^q(\rr^2)\cap L_y^qL^r_x(\rr^2)$, for any $1\leq r\leq\infty$ satisfying
\[
\frac{1}{r}+\frac{1}{2q}>\frac{1}{2}.
\]
\end{corollary}

Here is the main result about the spatial decay of the solitary-wave solutions.
\begin{theorem}
Let  $\sigma_0>0$ be in Lemma \ref{decay-expon-lem}. Then, for any $\sigma\in[0,\sigma_0)$
and any $0\leq s<3/2$, it transpires that $|x|^se^{\sigma|y|}\ff(x,y)\in
L^1\left(\rr^2\right)\cap L^\infty\left(\rr^2\right) $.
\end{theorem}
 \begin{proof}
  Without loss of generality, assume that
$s=0$. By using Lemma \ref{decay-expon-lem} and the proof of
Corollary 3.14 in \cite{bonali}, with natural modifications, it may be
seen that there is a $\ti{\sigma}\geq\sigma_0$ such that $e^{\sigma|y|}\ff(x,y)\in
L^1\left(\rr^2\right)$, for any $\sigma<\ti{\sigma}$.
The  inequality
\begin{equation}
|\ff(x,y)|e^{\sigma|y|}\leq\int_{\rr^2}|K(x-z,y-w)|e^{\sigma|y-w|}|\ff(z,w)|e^{\sigma|w|}|\ff(z,w)|^p\;dzdw
\end{equation}
and the facts $\ff(x,y)e^{\sigma|y|}\in L^1\left(\rr^2\right)$,
$\ff\in L^\infty\left(\rr^2\right)$ and $K(x,y)e^{\sigma|y|}\in
L^2\left(\rr^2\right)$, for any $\sigma<\sigma_0$, entails that
$\ff(x,y)e^{\sigma|y|}\in L^\infty\left(\rr^2\right)$, for the same range of
$\sigma$.
 \end{proof}
Finally, the following theorem deals with  the analyticity of the solitary-wave solutions.  Of
course, for this, one needs to restrict $p$ so that $z \mapsto z^p$ is analytic in a full
neighborhood of the origin in $\C$.
\begin{theorem}
Let $1\leq p<4$ be an integer. Then, there is a $\sigma>0$ and a holomorphic function $\mathpzc{f}$ of
two variables $z_1$ and $z_2$, defined in the domain
\[
\mathpzc{H}_\sigma=\left\{(z_1,z_2)\in\C^2\;;\;|\rm{Im}(z_1)|<\sigma,\;|\rm{Im}(z_2)|<\sigma\right\}
\]
such that $\mathpzc{f}(x,y)=\ff(x,y)$ for all $(x,y)\in\rr^2$.
\end{theorem}

Similar results are obtained by the same method for related evolution equations in \cite{BL} and
\cite{bonali}.    Results of this nature for dispersive equations
made via Gevrey-space analysis appear in  \cite{BG}
(and see also the reference therein).

\begin{proof}
 By the Cauchy-Schwarz inequality,  Theorem \ref{regularity} implies that $\what{\ff}\in L^1\left(\rr^2\right)$. Equation \eqref{2bozk} implies in turn that
\begin{gather}
|\xi|\left|\what{\ff}\right|(\xi,\eta)\leq
\overbrace{|\what{\ff}|\ast\cdots\ast |\what{\ff}|}^{p+1}(\xi,\eta),\\
|\eta|\left|\what{\ff}\right|(\xi,\eta)\leq \underbrace{|\what{\ff}|\ast \cdots\ast|\what{\ff}|}_{p+1}(\xi,\eta).
\end{gather}
Denote by $\mathscr{T}_1$ the  correspondence $\mathscr{T}_1(|\what{\ff}|)=|\what{\ff}|$ and, for
$m\geq1$,
$\mathscr{T}_{m+1}(|\what{\ff}|)=\mathscr{T}_m(|\what{\ff}|)\ast|\what{\ff}|$.
A straightforward induction yields
\begin{equation}
r^m|\what{\ff}|(\xi,\eta)\leq (m-1)!\;(p+1)^{m-1}\mathscr{T}_{mp+1}(|\what{\ff}|)(\xi,\eta),
\end{equation}
where $r=|(\xi,\eta)|$. It follows that
\[
\begin{split}
r^m|\what{\ff}|(\xi,\eta)&\leq(m-1)!\;(p+1)^{m-1}\left\|\mathscr{T}_{mp+1}(|\what{\ff}|)\right\|_{L^\infty(\rr^2)}\\
&\leq (m-1)!\;(p+1)^{m-1}\left\|\mathscr{T}_{mp}(|\what{\ff}|)\right\|_{L^2(\rr^2)}\|\what{\ff}\|_{L^2(\rr^2)}\\
&\leq (m-1)!\;(p+1)^{m-1}\|\what{\ff}\|_{L^1(\rr^2)}^{mp}\|\what{\ff}\|_{L^2(\rr^2)}^2.
\end{split}
\]
Let
\[
a_m=\frac{(p+1)^{m-1}\|
\what{\ff}\|_{L^1(\rr^2)}^{mp}\|\what{\ff}\|_{L^2(\rr^2)}^2}{m},
\]
so that
\[
\frac{a_{m+1}}{a_m}\longrightarrow(p+1)\|\what{\ff}\|_{L^1(\rr^2)}^p,
\]
as $m\to+\infty$.  In consequence, the series $\sum_{m=0}^\infty t^m
r^m|\what{\ff}|(\xi,\eta)/m!$ converges uniformly in
$L^\infty(\rr^2)$ provided
$0<t<\sigma=\frac{1}{p+1}\|\what{\ff}\|_{L^1(\rr^2)}^{-p}$. Hence
$e^{tr}\what{\ff}(\xi,\eta)\in L^\infty(\rr^2)$, for
$t<\sigma$. Now define the function
\[
\mathpzc{f}(z_1,z_2)=\int_{\rr^2}e^{i(\xi z_1+\eta z_2)}\what{\ff}(\xi,\eta)\,d\xi d\eta.
\]
By the Paley-Wiener Theorem, $\mathpzc{f}$ is well defined and
analytic in $\mathpzc{H}_\sigma$ while Plancherel's Theorem assures that  $\mathpzc{f}(x,y)=\ff(x,y)$ for all $(x,y)\in\rr^2$. This proves the theorem.
\end{proof}

\section*{Acknowledgments}
The financial support of the research council of Damghan university  with the grant number 93/MATH/125/227 is acknowledged by the first author. AP was partially supported by
CNPq.  Both AP and JB were supported by FAPESP during part of this collaboration.  JB, who
also received support from the University of Illinois at Chicago,
thanks the Instituto de Matem\'{a}tica, Estat\'{\i}stica e Computa\c{c}\~{a}o Cient\'{\i}fica
 at UNICAMP and the Instituto Nacional de Matem\'atica Pura e Aplicada for hospitality during
the development of this work.

\end{document}